\documentstyle[12pt]{article}                  % Use this for Latex (2.09)
\input{amssymb.sty}
\makeatletter
\def\nothing#1{}
\newdimen\earraycolsep
\renewcommand{\thetable}{\arabic{table}}
\renewcommand{\thefigure}{\arabic{figure}}
\renewcommand{\title}[1]{%
  \vspace*{120\p@}%
  {\parindent \z@ \raggedright \reset@font
    \bfseries #1\par
    \nobreak
    \vskip 36\p@
  }}
\def\author#1{{\pretolerance=10000 \raggedright \advance \leftskip by 1in
\noindent #1 \vskip 1pc}}
\def\affiliation#1{{\advance\leftskip by 1in \noindent #1 \vskip -1pc}}
\def\refnote#1{{$^{\hbox{\scriptsize #1}}$}}
\def\affnote#1{\llap{$^{\hbox{\scriptsize #1}}$}}

\renewcommand\section{\@startsection{section}{1}{\z@}{2pc \@plus
      1ex minus .2ex}{1pc \@plus .2ex}{\reset@font
      \normalsize\bfseries\noindent
      {\addtocounter{section}{1}}\arabic{section}\
      {\setcounter{subsection}{0}
      \setcounter{subsubsection}{0}\setcounter{equation}{0}} }}
\renewcommand\subsection{\@startsection{subsection}{2}{\z@}{1pc \@plus 1ex
    minus.2ex}{1pc \@plus .2ex}
    {\reset@font\normalsize\bfseries
    \noindent{\addtocounter{subsection}{1}}%
    {\setcounter{subsubsection}{0}}\arabic{section}.\arabic{subsection}\ }}
\renewcommand\subsubsection{\@startsection{subsubsection}{3}{\parindent}
        {1pc \@plus 1ex minus.2ex}{-0.5em}{\reset@font\normalsize\bfseries%
        {\addtocounter{subsubsection}{1}} \hspace*{.6cm}
        \arabic{section}.\arabic{subsection}.\arabic{subsubsection}
        \hspace*{-7mm}}}

\def\AmS{{\protect\the\textfont2%
        A\kern-.1667em\lower.5ex\hbox{M}\kern-.125emS}}

\def\p@LaTeX{{\family{times}\series{m}\shape{n}\selectfont
L\kern-.36em\raise.3ex\hbox{\scriptsize A}\kern-.15em
T\kern-.1667em\lower.7ex\hbox{E}\kern-.125emX}}

\newlength{\colwidth}

\def\@oddhead{\hfil}
\def\@evenhead{\hfil}
\def\@oddfoot{{\bfseries\hfil\thepage}}
\def\@evenfoot{{\bfseries\thepage\hfil}}
\def\fnum@figure{\footnotesize\raggedright{\bfseries \figurename~\thefigure.}}
\def\fnum@table{\normalsize\raggedright{\bfseries \tablename~\thetable.}}
\long\def\@makecaption#1#2{\vskip 10\p@ {#1 #2\par}}
\long\def\@makefntext#1{\setbox0=\hbox{$\m@th^{\@thefnmark}$}\noindent
\hangindent=\wd0 \box0 #1}
\newbox\@atbox
\long\def\atable#1#2#3{\begin{table}[tbp]\centering\footnotesize
\setbox\@atbox\hbox{#2}
\parbox{\wd\@atbox}{\caption{#1}}\par\smallskip
#2
\par\smallskip\parbox{\wd\@atbox}{\raggedright #3}
\end{table}}

\newtheorem{theorem}{Theorem}

\newtheorem{proposition}[theorem]{Proposition}
\newtheorem{lemma}[theorem]{Lemma}
\newtheorem{corollary}[theorem]{Corollary}

\newtheorem{definition}[theorem]{Definition}

\newtheorem{remark}[theorem]{Remark}

\newcommand{\ra}{\rightarrow}
\newcommand{\fl}{\forall}
\newcommand{\wt}{\widetilde}

\newcommand{\s}{\sigma}
\newcommand{\D}{\Delta}

\newcommand{\Zb}{\mathbb{Z}}
\newcommand{\Uc}{\mathcal{U}}
\newcommand{\ot}{\otimes}

\newcommand{\Hc}{\mathcal{H}}
\newcommand{\g}{\gamma}
\newcommand{\vp}{\varphi}
\newcommand{\ve}{\varepsilon}
\newcommand{\Cb}{\mathbb{C}}

\newcommand{\Fa}{\mathfrak{a}}

\newcommand{\Fg}{\mathfrak{g}}

\newcommand{\Fl}{\mathfrak{l}}

\def\Cb{{\mathbb C}}
\def\Gb{{\mathbb G}}

\def\Rb{{\mathbb R}}

\def\Zb{{\mathbb Z}}

\def\Ac{{\mathcal A}}
\def\Bc{{\mathcal B}}
\def\Cc{{\mathcal C}}
\def\Dc{{\mathcal D}}

\def\Gc{{\mathcal G}}
\def\Hc{{\mathcal H}}
\def\Ic{{\mathcal I}}
\def\Ic{{\mathcal I}}
\def\Jc{{\mathcal J}}
\def\Kc{{\mathcal K}}
\def\Lc{{\mathcal L}}
\def\Mc{{\mathcal M}}
\def\Nc{{\mathcal N}}
\def\Pc{{\mathcal P}}
\def\Qc{{\mathcal Q}}
\def\Rc{{\mathcal R}}

\def\Tc{{\mathcal T}}
\def\Uc{{\mathcal U}}
\def\Vc{{\mathcal V}}

\def\a{\alpha}
\def\b{\beta}
\def\c{\chi}
\def\d{\delta}

\def\g{\gamma}
\def\k{\kappa}
\def\i{\iota}
\def\lb{\lambda}
\def\om{\omega}
\def\s{\sigma}
\def\t{\theta}
\def\ve{\varepsilon}
\def\ph{\phi}
\def\vp{\varphi}

\def\vf{\varpi}
\def\z{\zeta}

\def\D{\Delta}
\def\G{\Gamma}
\def\Lb{\Lambda}
\def\Om{\Omega}

\def\Ph{\Phi}
\def\Ps{\Psi}
\def\Th{\Theta}

\def\fl{\forall}
\def\ify{\infty}
\def\lgl{\langle}
\def\nb{\nabla}
\def\op{\oplus}
\def\ot{\otimes}
\def\part{\partial}
\def\rgl{\rangle}
\def\sbs{\subset}
\def\semi{\rtimes}
\def\sm{\simeq}
\def\ts{\times}
\def\wdg{\wedge}

\def\ra{\rightarrow}
\def\longra{\longrightarrow}

\def\text{\hbox}

\def\Ker{\mathop{\rm Ker}\nolimits}

\def\Sign{\mathop{\rm Sign}\nolimits}

\def\Trace{\mathop{\rm Trace}\nolimits}

\def\boxit#1#2{\setbox1=\hbox{\kern#1{#2}\kern#1}%
\dimen1=\ht1 \advance\dimen1 by #1 \dimen2=\dp1 \advance\dimen2 by #1
\setbox1=\hbox{\vrule height\dimen1 depth\dimen2\box1\vrule}%
\setbox1=\vbox{\hrule\box1\hrule}%
\advance\dimen1 by .4pt \ht1=\dimen1
\advance\dimen2 by .4pt \dp1=\dimen2 \box1\relax}

\def\build#1_#2^#3{\mathrel{
\mathop{\kern 0pt#1}\limits_{#2}^{#3}}}

\def\@nbibitem#1{\noindent \hangindent=2pc \hangafter=1
\refstepcounter{enumi}\hbox to 2pc{\arabic{enumi}.\hfil}%
\immediate\write\@auxout{\string\bibcite{#1}{\arabic{enumi}}}}
\def\numbibliography{%
\section*{REFERENCES}%
\bgroup\footnotesize
\setcounter{enumi}{0}%
\def\newblock{\hskip .11em plus.33em minus.07em}%
\let\bibitem\@nbibitem}
\def\endnumbibliography{\par\egroup}
\makeatother
\begin{document}

\begin{center}
{ \sc Differentiable cyclic cohomology and Hopf algebraic structures
                    in transverse geometry}
\end{center}
\medskip
\author{\bf Alain Connes\refnote{1} and 
    Henri Moscovici\refnote{2}\footnote{The second author has 
    conducted part of this work during employment by the Clay Mathematics 
    Institute, as a CMI Scholar in residence at Harvard University. 
    His research is supported in part by the U.S. National Science Foundation
    award no. DMS-9988487 and the U.S.-Israel Binational Science Foundation
    award no. 1997405.
 }}

\affiliation{\affnote{1}  Coll\`ege de France,
3, rue Ulm,
75005 Paris\\
and\\
I.H.E.S., 35, route de Chartres, 91440 Bures-sur-Yvette \\
\affnote{2} Department of Mathematics,
The Ohio State University \\
231 West 18th Avenue, Columbus, OH 43210-1174, USA  
}

\vspace{1cm}

\begin{center}
{\it Dedicated to Andr\'e HAEFLIGER, whose {\em philosophy} 
of differentiable cohomology has inspired the present paper}
\end{center}

\medskip

\begin{abstract}
We prove a cyclic cohomological analogue of Haefliger's
van Est-type theorem for the groupoid of germs of diffeomorphisms
of a manifold. The differentiable version of cyclic cohomology is
associated to the algebra of transverse differential operators on 
that groupoid, which is shown to
carry an intrinsic Hopf algebraic structure.
We establish a canonical isomorphism
between the periodic Hopf cyclic cohomology of this extended Hopf
algebra 
and the Gelfand-Fuchs cohomology of the Lie algebra of formal vector
fields. We then show that this isomorphism can be explicitly
implemented at the cochain level, by a cochain map
constructed out of a fixed torsion-free linear connection. 
This allows the direct treatment of the index formula for
the hypoelliptic
signature operator -- representing the diffeomorphism invariant
transverse fundamental $K$-homology class of an oriented manifold 
-- in the general case, when this operator is
constructed by means of an arbitrary coupling connection.
\end{abstract}

\newpage

\centerline{\bf Introduction}
\medskip

The local index formula for hypoelliptic differential
operators in a diffeomorphism invariant setting \cite{CM1} 
expresses 
the Chern character of such operators in terms of a particular
kind of cocycles of the cyclic bicomplex, satisfying a property
analogous to that involved in Haefliger's definition of
differentiable cohomology for the groupoid of germs of 
diffeomorphisms on a manifold \cite{H2}. 
Applying it to the transverse index problem on foliations,
we have shown in \cite{CM2} that, if the
complete transversal is chosen to be an affine flat manifold, 
the ``differentiable'' cocycle representing 
the Chern character of such an operator is automatically
in the range of a natural characteristic map, originating
from the cyclic cohomology of a Hopf algebra $\Hc_{n}$ canonically
associated to $ {\rm Diff} \ (\Rb^{n}  ) $. In turn,
the periodic Hopf cyclic cohomology of $\Hc_{n}$ 
and its $SO(n)$-relative version
were shown to be canonically isomorphic to corresponding
Gelfand-Fuchs cohomologies of formal vector fields. The
upshot was proving that, in cohomological form, the
index formula for transversely hypoelliptic operators on foliations
can be expressed in terms of Gelfand-Fuchs classes, transported
via the characteristic map.
In this paper we incorporate ab initio the curvature into our 
approach and thus dispense with the geometrically unsatisfactory
flatness condition on the transversal.

The new device that allows the direct treatment of the curved case 
is a ``thickened'' version ${\Hc}_{FM}$
of the Hopf algebra  $\Hc_{n}$. If one regards $\Hc_{n}$ as
an algebra of transverse differential operators with constant coefficients,
${\Hc}_{FM}$ should be viewed as the algebra of transverse differential
operators with variable coefficients on $FM \bar{\semi} \G_{M} \,$,
the \'etale groupoid of germs of diffeomorphisms of a 
given $n$-manifold $M$ lifted to its frame bundle $FM$.
The algebra ${\Hc}_{FM}$ is naturally a bimodule over the
coefficient ring $ \Rc_{FM} = C^{\ify} (FM) $ and it affords
a Hopf structure only in an extended sense, the coproduct
taking values in the tensor product over $ \Rc_{FM} $.
Like its precursor $\Hc_{n}$, ${\Hc}_{FM}$ too gives rise to a natural
cyclic module, consisting of $ (\Rc_{FM}€, \Rc_{FM}€)$-bimodules of
multidifferential operators.
 
The main result of this paper, which can be viewed as a
cyclic cohomological analogue of Haefliger's van Est-type
theorem for the groupoid of germs of diffeomorphisms of
$M \, $  \cite[Theorem IV.4]{H2},
establishes a canonical
isomorphism between the periodic Hopf cyclic cohomology of ${\Hc}_{FM}$ 
and the Gelfand-Fuchs cohomology of the Lie algebra of formal vector
fields on $\Rb^n$, as well as between their corresponding
relative versions. This isomorphism is implemented  at the cochain 
level by a cochain map manufactured out of a given
torsion-free connection. In this way we exhibit
a purely geometric construction, tracking the displacement of
the connection and curvature forms
under diffeomorphisms, for the cyclic cohomological
counterparts of the classical
Chern and secondary classes of foliations (\cite{Bo}, \cite{H1}).
\bigskip

\noindent {\it Acknowledgment.} We express our thanks to Georges
Skandalis, for his critical reading of the manuscript and for his
insightful comments.

\section{The algebra ${\Hc}_{FM}$}

We begin by introducing ${\Hc}_{FM}$, as an algebra of
transverse differential operators on the smooth \'etale groupoid
associated to the diffeomorphisms of the base manifold
acting on its frame bundle.

Let $M$ be a smooth $n$-manifold and let $FM$ 
denote the bundle of frames on $M$.
We denote by $\G_{M}$ the pseudogroup of all
local diffeomorphisms of $M$; its elements are
partial diffeomorphisms 
$ \psi : {\rm Dom} \psi \ra {\rm Ran} \psi $,
with both the domain ${\rm Dom} \psi$ and range ${\rm Ran} \psi$
open subsets of $M$. 
The prolongation
to $FM$ of a partial diffeomorphism $\psi \in \G_{M}€$
will be denoted
by $ \wt{\psi}$. The set of all pairs $ (u, \wt{\vp}) $ with
$\vp \in \G_{M}$ and $u \in {\rm Ran}  \wt{\vp} $ will be denoted
$FM \semi \G_{M} $. 
We next form the associated smooth \'etale groupoid
of germs of lifted diffeomorphisms
$$
FM \bar{\semi} \G_{M}€: = \{ \, [u, \wt{\vp}] \, ; \quad \vp \in \G_{M}, 
\quad u \in {\rm Ran}  \wt{\vp} \,  \} ,
$$
where $ [u, \wt{\vp}] $ stands for the class of $(u, \wt{\vp}) \in
FM \semi \G_{M} $ with respect to the equivalence relation which
identifies $(u, \wt{\vp})$ and $(v, \wt{\psi})$ if $ u = v$ and
$ {\wt{\vp}}^{-1} $ coincides with ${\wt{\psi}}^{-1}$ on a neighborhood
of $u$. The corresponding
{\it source} and {\it target} maps are 
$$
s [u, \wt{\vp}]  = {\wt{\vp}}^{-1}€(u) \in {\rm Dom} \wt{\vp} \, ,
\quad \hbox{resp.} \quad
t [u, \wt{\vp}] = u \, ,
$$
and the composition rule is
$$
[u, \wt{\vp}] \circ [v, \wt{\psi}] = [u, \wt{\vp} \circ \wt{\psi}] \, ,
\quad \hbox{if} \quad v \in {\rm Dom} \wt{\vp} \quad \hbox{and}
\quad \wt{\vp} (v) = u \, . 
$$
We let 
$$
   \Ac \equiv {\Ac}_{FM}€ : = C_c^{\ify} (FM  \bar{\semi} \G_{M}€ )
$$
denote its convolution algebra. From a practical standpoint,
it is convenient to regard ${\Ac}_{FM}$ as being linearly spanned
by monomials of the form
$$
f \, U_{\psi}^* \, , \quad \hbox{with} \quad
f \in C_c^{\ify} ({\rm Dom} \wt{\psi}) \, ,
$$
where the asterisk stands for the inverse, with the understanding that
$$ f_1 \, U_{\psi_1}^* \, \equiv f_2 \, U_{\psi_2}^* \, \quad
\hbox{iff} \quad f_1 = f_2 \quad \hbox{and} \quad 
\wt{\psi}_1 \vert_V= \wt{\psi}_2 \vert_V \, , 
$$
where $V$ is a neighborhood of ${\rm Supp} (f_i)$, $i=1$ or $2$.
The multiplication rule for such monomials is given by
$$
f_1 \, U_{\psi_1}^* \cdot f_2 \, U_{\psi_2}^* = 
f_1 (f_2 \circ \wt{\psi}_1) \,
U_{\psi_2 \psi_1}^* \, ;
$$
note that, by hypothesis, the support of $f_1 (f_2 \circ \wt{\psi}_1)$ 
is a compact subset of 
${\rm Dom}  \wt{\psi_1} \cap {\wt{\psi}_1}^{-1} \, ({\rm Dom} \wt{\psi}_2 )
\sbs {\rm Dom} ( \wt{\psi}_2 \circ \wt{\psi}_1 ) \, $.

The function algebra 
$$ \Rc \equiv {\Rc}_{FM} : = C^{\infty}€(FM)
$$ 
acts in two ways on $\Ac $, 
by left multiplication operators
\begin{equation} \label{lb}
\a (b) (f \, U_{\psi}^* ) = b \cdot f \, U_{\psi}^* \, , \quad
b \in \Rc \, , 
\end{equation}
and by right multiplication operators
\begin{equation} \label{rb}
\b (b) (f \, U_{\psi}^* ) =  \ f \, U_{\psi}^* \cdot b 
= b \circ \wt{\psi} \cdot f \, U_{\psi}^* \, , \quad
b \in \Rc \, . 
\end{equation}
On the other hand, any vector field $Z$ on $F$ can be extended to
a linear transformation, although no longer a derivation in general,
$Z \in \Lc \, ({\Ac})$, by setting 
\begin{equation}  \label{Z}
Z (f \, U_{\psi}^* ) = Z (f) \, U_{\psi}^* \, , \quad 
f \, U_{\psi}^* \in {\Ac} .
\end{equation}

The following definition can be easily adapted to cover the case of
an arbitrary base manifold. However, for the simplicity of the
exposition,
we shall {\it assume from now on that the manifold $M$ admits a
finite atlas}.

\begin{definition} A {\em transverse differential operator} on
the groupoid $FM  \bar{\semi} \G_{M}$ is an element of the subalgebra 
of linear operators on $\Ac$
$$ \Hc \equiv {\Hc}_{FM}€ \sbs \Lc€ \, ({\Ac}_{FM}€)
$$
generated by the 
transformations (\ref{lb}), (\ref{rb}) and (\ref{Z}). 
A {\em $p$-differential operator} on $FM  \bar{\semi} \G_{M}$
is a $p$-linear transformation $H$ on ${\Ac}_{FM}$ with values 
in ${\Ac}_{FM}$, of the form
\begin{equation} \label{mdiff} 
	H( a^1 ,..., a^p ) = \sum_{i=1}^r h_i^1 (a^1 ) \cdots h_i^p ( a^p )
	\, , \quad a^1 ,..., a^p \, \in {\Ac}_{FM} \, ,
\end{equation}
with $h_i^1, \ldots, h_i^p \, \in {\Hc}_{FM}$.
\end{definition}
\noindent The adjective {\it transverse} is meant to emphasize the
distinction between this notion and that of {\it longitudinal} or
{\it invariant} (pseudo)differential operator of \cite{C0} 
(cf. also \cite{N-W-X} 
for a systematic treatment of the latter in the context of Lie
algebroids.) 
\smallskip

 There are two built-in algebra homomorphisms $\a : \Rc \ra \Hc$ 
and  $  \b : \Rc \ra \Hc$, whose images commute
\begin{equation} \label{ab}
 \a (b_1) \ \b (b_2) = \b (b_2) \ \a (b_1) \, , \quad
	\fl \quad  b_1, b_2 \in \Rc ; 
\end{equation}
we shall view $\Hc$ as an $(\Rc, \Rc)$-bimodule with the
{\it left action} of $\Rc$ given by {\it left multiplication via} $\a$
and the {\it right action} given by {\it left multiplication via} $\b$.
More generally, for any $p \geq 1$,
the space of  $p$-differential operators on $FM  \bar{\semi} \G_{M}$
\begin{equation} \label{hp}
{\Hc}^{\, [p]}€ \equiv {\Hc}_{FM}€^{\, [p]}€
\end{equation}
can be endowed with an
$(\Rc, \Rc)$-bimodule structure in a similar fashion, via
left multiplication by means of $\a$, respectively $\b$.

The remainder of this section will be devoted to the proof
of two pivotal results,
clarifying the structure of these bimodules.  

We fix a torsion-free connection on FM, with connection form
$\om = ( \om^i_j )$. We denote by $\t = ( \t^i )$ the canonical
form
$$ \t_u (Z) = u^{-1} ( \pi_{*} (Z) ) \, , \quad \fl Z \in T_u FM \, ,
$$
where $\pi : FM \ra M$ is the projection and the frame $u \in FM$ is
viewed as an isomorphism of vector spaces
$u : \Rb^n \build \longra_{}^{\simeq} T_{\pi(u)} M $.

Let $X_1 , ..., X_n $ be the standard horizontal vector fields
corresponding to the standard basis of $\Rb^n$ and 
$\{ Y_i^j \}$ the fundamental vertical vector fields  
corresponding to the standard basis of ${\Fg \Fl} (n, \Rb )$, such that
$\{ X_k , Y_i^j \} $ and $\{ \t^k , \om^i_j \}$ are dual to
each other. Since $\{ X_k \vert_u , Y_i^j \vert_u  \} $ form a basis
of $T_u FM $ for every $u \in FM$, as an algebra, $\Hc$ is generated
by $\a ( \Rc )$, $\b ( \Rc )$ and the endomorphisms of $\Ac$
associated to the vector fields $X_k$ and $Y_i^j$ (cf. (\ref{Z})).
The usual commutation relations for the vector fields associated
to a torsion-free connection continue to hold in $\Hc$:
\begin{equation} \label{YX}
\matrix{
\displaystyle  [Y_i^j , Y_k^{\ell}] = \d_k^j Y_i^{\ell} - \d_i^{\ell} Y_k^j 
\, , \hfill \cr \cr
\displaystyle  [Y_i^j , X_k] = \d_k^j X_i \, , \hfill \hfill \cr \cr
\displaystyle  [X_k , X_{\ell}] = \sum \ \a (R_{jk\ell}^i ) Y_i^j \, ,
\hfill \cr
}  
\end{equation}
where the functions $R_{jk\ell}^i \in \Rc$ are related to the curvature
form $\Om$ of the given connection by
$$ \Om_j^i = \sum_{k < \ell} \ R_{jk\ell}^i \ \t^k \wedge \t^{\ell} \, ,
 \quad R_{jk\ell}^i = - R_{j\ell k}^i \, .
$$
Also, for any $b \in \Rc$ one has
\begin{equation} \label{leftb}
\matrix{
\displaystyle [Y_i^j , \a (b)] &= &\a (Y_i^j \, (b)) \, , \hfill \cr \cr
\displaystyle [X_k , \ \a (b)] &= &\a (X_k \, (b)) \ . \hfill \cr 
} 
\end{equation}
On the other hand, the commutators between the
standard horizontal vector fields
and $\b (\Rc)$ introduce new operators acting on $\Ac$. 
Indeed, while the 
fundamental vector fields are $\G_{M}€$-invariant
and therefore
\begin{equation} \label{Yb} 
   [Y_i^j , \b (b)] = \b (Y_i^j \, (b))  \, , 
\end{equation}
for the standard horizontal vector fields one has, with the
usual summation convention,
\begin{equation} \label{Xd}
 U_{\vp} \, X_k \, U_{\vp}^* - X_{k}€ =  \ \rho_{jk}^i \, Y_i^j \ ;
\end{equation}
the coefficients $\rho_{jk}^i \, $, which are functions
on $FM \bar{\semi} \G_{M}€$ involving the 
second order jet of the local diffeomorphism,
can be alternatively expressed as 
\begin{equation} \label{g}
 \rho_{jk}^i \, = \, \g_{jk}^i \circ \wt{\vp}^{-1} \, , 
 \quad \hbox{with} \quad
 \g_{jk}^i  =  \langle \wt{\vp}^{*}€\om_{j}€^{i}€ \, ,
 X_{k}€\rangle \, .
\end{equation}
It then follows that
\begin{equation} \label{Xb}   
  [X_{k} , \b (b)] = \b (X_{k} (b)) + \ \b (Y_i^j (b)) \, \d_{j k}^i \ , 
\end{equation}
where $\d_{jk}^i \in \Lc \, ({\Ac}) $ are the operators of multiplication
by the functions $\g_{jk}^i \in C^{\ify} (FM \bar{\semi} \G_{M})$, 
\begin{equation} \label{d} 
\d_{jk}^i (f \, U_{\vp}^*) := \g_{jk}^i \cdot f \, U_{\vp}^* \,
 = f \, U_{\vp}^* \cdot \rho_{jk}^i \, . 
\end{equation}

\smallskip

\begin{lemma} \label{dijk}
The operators $\d_{jk}^i $ are transverse differential operators on
the \'etale groupoid $FM  \bar{\semi} \G_{M}$.
\end{lemma}

\noindent{\it Proof.} Due to the existence of a finite atlas on $M$ 
and using a partition of unity argument, it suffices to show that
$$\b (\chi \circ \pi) \ \d_{jk}^i \in \Hc$$
 for any function
 $\chi \in C_c^{\ify} (M)$ with support in a coordinate chart
 $U \sbs M$. Equation (\ref{Xb}) ensures that
\begin{equation} \label{L1}
 \b \left( Y_i^j (b) \right) \, 
 \b (\chi \circ \pi) \, \d_{j k}^i 
 \in \Hc \, , \quad \fl b \in \Rc \, .
\end{equation}
Choose $b \in M_n ( \Rc ) $ such that,  
in local coordinates $ u = (x^k, \, y_i^j ) \, $ 
on $\pi^{-1} (U)$,
$$
b (u) = y \, .
$$
Then
\begin{equation} \label{L2}
\matrix{
Y_i^j (b_{\ell}^s) \, \mid_{\ \pi^{-1} (U)} \,
=  \ y_i^r \ \frac{\partial}{\partial 
\, y_j^r} \ (y_{\ell}^s) =  y_i^s \, \d_{\ell}^j \, . \hfill \cr
}
\end{equation}
We now choose $d \in M_n ( \Rc )$ such that its
restriction to  $\pi^{-1} (U)$ is
$$
	d (u) =  y^{-1} .
$$
From (\ref{L2}) it follows that
$$
\b ( d_s^r ) \,  \,  \b \left ( \ Y_i^j (b_{\ell}^s ) \right) \, 
\b (\chi \circ \pi) \, \d_{jk}^i \, = \, 
\b (\chi \circ \pi) \, \d_{lk}^r \, ,
$$
and by  (\ref{L1}) the left hand side
belongs to the algebra $\Hc$. $ \qquad \Box$
\medskip

An equivalent form for the equation (\ref{g}) is
\begin{equation} \label{g'}
 \wt{\vp}^{*} \om \, - \, \om \, = \, \g \cdot \t \, ;
\end{equation} 
from this it easily follows that,
as a vector valued function on the frame bundle, $\g$ is 
tensorial with respect to the right action of the
structure group $GL^+ (n, \Rb)$. At the infinitesimal level,
this translates into the following expressions
for the commutators between the 
$Y_j^i$ and $\d_{\ell m}^k$:
\begin{equation} \label{Yd}
\displaystyle [Y_j^i , \d_{\ell m}^k ] = \d_{\ell}^i \ \d_{jm}^k
\, + \, \d_m^i \ \d_{\ell j}^k \, - \, \d_j^k \ \d_{\ell m}^i \, ,
\end{equation}
The commutators with $X_k$ however yield new
operators, involving higher 
order jets of diffeomorphisms:
\begin{equation} \label{highd}
\d_{jk,\ell_1 \ldots , \ell_r}^i \, = \,
[X_{\ell_r} , \ldots [X_{\ell_1} , \d_{jk}^i] \ldots ] \,;
\end{equation}
these operators, acting on $\Ac$, have the form
\begin{equation} \label{d'} 
\d_{jk,\ell_1 \ldots , \ell_r}^i \, ( f \, U_{\vp}^*) \,= \,
 \g_{jk, \ell_1 \ldots \ell_r}^i \cdot f \, U_{\vp}^* \, ,
\qquad  \g_{jk, \ell_1 \ldots \ell_r}^i \, =
\, X_{\ell_r} \cdots X_{\ell_1} (\g_{jk}^i). 
\end{equation}
In particular, they form an {\it abelian subalgebra} and also
commute with image of $\Rc$ through both maps $\a$ and
$\b$.
\smallskip

We are now ready to prove the main results of this section.
Recall that $\Hc$ is viewed as an $(\Rc , \Rc)$-bimodule under
the {\it left action} of $\Rc$ via multiplication by the
image of $\a$, resp. $\b$. 

\begin{proposition} \label{free}
The $(\Rc , \Rc)$-bimodule $\Hc$ is free over $\Rc \ot \Rc$.
The choice of a torsion-free connection on $FM$ gives rise
to a Poincar\'e-Birkhoff-Witt-type basis of $\Hc$ over
$\Rc \ot \Rc$.
\end{proposition}

\noindent {\it Proof}. Once the connection $\om$ is fixed, we use the
same notation as above for the associated generators of the algebra $\Hc$. 
In addition, we shall need to employ two kinds of multi-indices,
whose entries form an increasing sequence with respect to
the obvious lexicographic order. The first kind are of the form
$$
I \, = \, \left\{ i_1 \leq \cdots \leq i_p \, ; \ 
\left( {j_1 \atop k_1}
\right) \leq \cdots \leq \left( {j_q \atop k_q} \right) \right\} \, ,
$$
while the second kind of the form
$$ \k \, = \,
\left\{ \left( {i_1 \hfill \atop j_1 \, k_1 ; 
\ell_1^1 \leq \ldots \leq \ell_{p_1}^1 } 
\right) \leq \cdots \leq \left( {i_r \hfill \atop j_r \, k_r ; 
\ell_1^r \leq \ldots 
\leq \ell_{p_r}^r} \right) \right\} \, .
$$
With this notation, we set
$$
Z_I = X_{i_1} \ldots X_{i_p} \, Y_{k_1}^{j_1} \ldots Y_{k_q}^{j_q} 
\quad \hbox{and} \quad
\d_{\k} = \d_{j_1 \, k_1 ; \ell_1^1 \ldots \ell_{p_1}^1}^{i_1} \ldots \d_{j_r 
\, k_r ; \ell_1^r \ldots \ell_{p_r}^r}^{i_r} \, .
$$

Using the relations (\ref{YX})-(\ref{Yb}), (\ref{Xb}), 
Lemma~\ref{dijk}, (\ref{Yd}) and (\ref{highd}), it is easy to check that
the collection $\{ \d_{\k} \cdot Z_I  \} \, $, where $I$ and $\k$ 
are multi-indices
of the first and second kind respectively, forms a generating set of $\Hc$
over $\Rc \ot \Rc$. 
We thus only need to prove that
\begin{equation} \label{eq2.1}
\sum_{I,\k} \ \a (\ell_{I,\k}) \, \b (r_{I,\k}) \, \d_{\k} \cdot Z_I =  0 
\quad \Rightarrow \quad
\ell_{I,\k} \ot  r_{I,\k} = 0  \, , \quad \fl \, (I, \k) \, .
\end{equation} 
Evaluating the expression in (\ref{eq2.1}) on an arbitrary
element $ \, f \, U_{\vp}^* \in \Ac \,$ one gets for any
$ \,  u \in FM$
\begin{equation} \label{eq2.2}    
\sum_{I,\k} \, \ell_{I,\k} (u) \, r_{I,\k} (v) \, \g_{\k} (u, 
 \wt{\vp} ) \, (Z_I \, f) (u) = 0 \, , \quad v = \wt{\vp} (u) \, .
\end{equation} 

Let us fix, for the moment, $u_0 , v_0 \in FM$ and set 
$$ \G^{(1)}€ \, (u_0 , v_0) \, 
= \, \{ \vp \in {\G}_{M}€ \, ; \ \wt{\vp} (u_0) = v_0 \} \, . 
$$
By varying $\vp \in \G^{(1)}€ \, (u_0 , v_0) $ ,
we shall first prove that (\ref{eq2.2}) implies that, 
for any fixed multi-index of the second kind $\k \,$  and
for any function $f \in C_{c}€^{\ify}€ (FM)$ supported around $u_{0}$,
one has
\begin{equation} \label{eq2.3}    
\sum_{I} \, \ell_{I,\k} (u_0) \, r_{I,\k} (v_0) \,
 (Z_I \, f) (u_0) = 0 \, .
\end{equation} 
Let $\Phi_{u_0} : \Rb^n \build \longra_{}^{u_0} T_{x_0} \, M \build 
\longra_{}^{\exp_{x_0}} M$ be normal coordinates around $x_0 = \pi (u_0)$ and 
similarly around $y_0 = \pi (v_0)$. Clearly,
$\vp \in \G^{(1)}€ \, (u_0 , v_0) $ iff
$\Phi_{v_0}^{-1} \circ \vp 
\circ \Phi_{u_0} \in {\rm Diff}_0^{(1)} :=$ the pseudogroup of all
local diffeomorphisms $\psi$ of $\Rb^{n}€$
such that their first order jet at $0$
satisfies $\ j_0^1 (\psi) = j_0^1 ({\rm Id}) \} \,$ .

In normal coordinates and using the customary notation
(see \cite{CM2}, also \cite{Wu}, for more detailed local 
calculations), 
one has:
\begin{equation} \label{eq2.y}
Y_i^j = y_i^{\mu} \ \frac{\partial}{\partial \, y_j^{\mu}} \equiv y_i^{\mu} \, 
\partial_{\mu}^j \, ,
\end{equation}
\begin{equation} \label{eq2.4}
X_k = y_k^{\mu} \, (\partial_{\mu} - \Gamma_{\a \mu}^{\nu} \, y_j^{\a} \, 
\partial_{\nu}^j ) \, ,
\end{equation} 
and
\begin{equation} \label{eq2.5}
\g_{ij}^k (x, y, \wt{\psi}) = \left( \wt{\G}_{\a\mu}^{\nu} (x) - 
\G_{\a \mu}^{\nu} (x) 
\right) \, y_j^{\a} \, y_i^{\mu} \, (y^{-1})_{\nu}^k
\end{equation} 
where
\begin{equation} \label{eq2.6}
\matrix{
\displaystyle \wt{\G}_{\a\mu}^{\nu} (x) &= &(\partial 
\psi (x)^{-1})_{\d}^{\nu} \,
\G_{\ve 
\z}^{\d} (\psi (x)) \, \partial_{\a} \, \psi (x)^{\ve} \, 
\partial_{\mu} \, \psi 
(x)^{\z} \cr \cr
&+ & \displaystyle (\partial \, \psi (x)^{-1})_{\d}^{\nu} \, 
\partial_{\mu} \, 
\partial_{\a} \, \psi (x)^{\d} \, . \hfill \cr
}
\end{equation}
Since, by the very choice of local coordinates,
$\G_{\a  \mu}^{\nu} (0) = 0 \, $, it follows from 
(\ref{eq2.6}) that
$$
 \wt{\G}_{\a\mu}^{\nu} (0) = \d_{\d}^{\g} \, \G_{\ve \z}^{\d} (0) \, 
\d_{\a}^{\ve} \, \d_{\mu}^{\z} + \d_{\d}^{\g} (\partial_{\mu} 
\, \partial_{\a} 
\, \psi^{d}) (0) \,
= \,  (\partial_{\mu} \, \partial_{\a} \, \psi^{\nu}) (0) \, .
$$
Therefore, by (\ref{eq2.5}),  $\fl \, \psi \in {\rm Diff}_0^{(1)} 
\, $, one has 
\begin{equation} \label{eq2.7}
\g_{ij}^k (0, y, \wt{\psi}) = \partial_{\mu} \, \partial_{\a} \, 
\psi^{\nu} (0) \cdot 
y_j^{\a} \, y_i^{\mu} \, (y^{-1})_{\nu}^k \, ,
\end{equation}
just like in the affine flat case \cite{CM2} .

We can now prove (\ref{eq2.3}) by induction on the ``height'' $\vert 
\k \vert$ 
of the multi-index $\k$, i.e. the number $p$ of indices in $\d_{jk ; \ell _{1}
 \ldots  \ell_p}^i $, counting the number of commutators with 
$X_{\ell}€$.

From (\ref{eq2.4})-(\ref{eq2.6}) it follows that the top degree 
component of the jet of $\psi$ occuring in the expression
$$
X_{\ell_p} \ldots X_{\ell_1} \, (\g_{jk}^i) (0, y, \wt{\psi})
$$
is of the form
$$
(\partial_{\ell_1} \ldots \partial_{\ell_p} \, \psi^{\nu} (0)) \cdot 
( y_{\cdot}^{\cdot}€€ ) \cdots 
( y^{-1} )_{\cdot}^{\cdot}€€  \, .
$$
By first choosing $\psi$ with $j_0^{p-1} (\psi) = j_0^{p-1} ({\rm Id})$, 
we can reduce to the situation $\vert \kappa \vert = p$, when
the equations look just like in the flat case. We can then
vary the $p$-th jet of $\psi$ to get the vanishing
of its coefficient and thus drop the height. 

We are now reduced to
proving that if $\fl \, u , v \in FM$,
$$
\sum_I \ \ell_I (u) \, r_I (v) \, (Z_I \, f) (u) = 0 \, , \ \fl \, f \in 
C_c^{\ify} (FM) \, ,
$$
then, for all indices $I$ of the first kind, one has
$$
\ell_I \ot r_I \, = \, 0 \, .
$$
This can be done by induction, in the same manner as before, using
(\ref{eq2.y}) and (\ref{eq2.4}). Alternatively, it also follows from
the Poincar\'e-Birkhoff-Witt theorem 
for Lie algebroids (\cite{Ri}, \cite{N-W-X}), applied to $TFM$.
 $ \qquad \Box$
\medskip

\begin{proposition} \label{multid}
There is a unique isomorphism of vector spaces
$$ T: \Hc^{ \, \ot_{\Rc}€ \, p}€ \, \equiv \, 
\underbrace{\Hc \ot_{\Rc}€ \ldots 
\ot_{\Rc}€ \Hc}_{p\hbox{-times}} \quad 
\build \longra_{}^{\simeq} \quad {\Hc}^{\, [p]}€ \, ,
$$
where ${\Hc}^{\, [p]}€$ (cf. (\ref{mdiff}), (\ref{hp}))
is the space of transverse $p$-differential operators,
such that $ \, \fl \, h_{1}€, \ldots , h_{p}€ \in \Hc \,$ and
$\, \fl \, a^{1}, \ldots , a^{p} \in \Ac \,$,
\begin{equation} \label{T}
T(h_{1}€ \ot_{\Rc} \ldots \ot_{\Rc} h_{p}€) \, (a^{1}€, \ldots , a^{p}€) 
\, = \,
 \  h_{1}€(a^{1}€) \cdots h_{p}€(a^{p}) \, .
\end{equation}
\end{proposition}

\noindent{\it Proof}. The assignment given by (\ref{T}) clearly 
extends by linearization to a well-defined epimorphism 
$ T: \Hc^{ \, \ot_{\Rc}€ \, p}€ \ra {\Hc}^{\, [p]}€ \,$. 
It remains to prove that $ {\rm Ker} \ T \, = \, 0 \, $.
Assume therefore that
\begin{equation} \label{nul}
\sum_{i}€ \, h_{1}€^{i}€(a^{1}€) \cdots h_{p}€^{i}€(a^{p}€) 
\, = \, 0 \, , \quad \fl \, a^{1}, \ldots , a^{p} \in \Ac \, .
\end{equation}
After fixing a Poincar\'e-Birkhoff-Witt $(\Rc , \Rc)$-basis  
$\, \{ B_{J} = \d_{\k} \cdot Z_I € \, ; J = I \cup \k \} \,$ of $\Hc $
as above,
we may express each $h_{j}€^{i}€$ under the form
$$
h_j^i = \sum_{J} \ \a (\ell_j^{i,J}) \, \b (r_j^{i,J}) \, 
B_{J} \, , \quad
\hbox{with} \quad \ell_j^{i,J} , \, r_j^{i,J} \in \Rc \, .
$$
By the hypothesis (\ref{nul}), for any
$ \, f_1 \, U_{\vp_1}^* , \ldots , f_p \, U_{\vp_p}^* \in 
\Ac$,
$$
\sum \ \ell_1^{i,J_1} \, B_{J_1} (f_1 \, U_{\vp_1}^*) \, r_1^{i,J_1} \, 
\ell_2^{i,J_2} \, B_{J_2} (f_2 \, U_{\vp_2}^*) \ldots r_{p-1}^{i,J_{p-1}} \, 
\ell_p^{i,J_p} \, B_{J_p} (f_p \, U_{\vp_p}^*) \, r_p^{i,J_p} \equiv 0 \, .
$$
Evaluating the left hand side at $u_0 \in FM$ and setting
$$
u_1 = \wt{\vp}_1 (u_0) \, , \ u_2 = \wt{\vp}_2 (u_1) \, , \ldots , u_p = 
\wt{\vp}_p (u_{p-1})
$$
one gets:
$$ \matrix{
\sum \ell_1^{i,J_1} (u_0) \, r_1^{i,J_1} (u_1) \, \ell_2^{i,J_2} (u_1) \ldots 
r_{p-1}^{i,J_{p-1}} (u_{p-1}) \, \ell_p^{i,J_p} (u_{p-1}) \, r_p^{i,J_p} (u_p)
\cr
\cdot \g_{\k_1} (\vp_1 , u_0) \, \g_{\k_2} (\vp_2 , u_1) \ldots 
\g_{\k_p} (\vp_p , u_{p-1}) \cr
\cdot Z_{I_1} (f_1) (u_0) \, Z_{I_2} (f_2) (u_1) \, \ldots \, Z_{I_p} (f_p) 
(u_{p-1}) \, = 0 \, .
}
$$
Following the same line of arguments as in the preceeding proof, 
we infer that $\fl \, J_1 , \ldots , J_p$ one has
\begin{equation} \label{nul*}
\sum_i  \ell_1^{i,J_1} (u_0) \, r_1^{i,J_1} (u_1) \, \ell_2^{i,J_2} (u_1) 
\ldots 
r_{p-1}^{i,J_{p-1}} (u_{p-1}) \, \ell_p^{i,J_p} (u_{p-1}) \, r_p^{i,J_p} 
(u_p)  = 0 .
\end{equation}

We now choose a basis $\{ b_{\lb} \}_{\lb = 1}^N$ of the finite-dimensional
subspace of $\Rc$ 
generated by the functions $\ell_k^{i,J_k} , r_k^{i,J_k} \, $ and
express them as linear combinations, with constant
coefficients, of the basis elements:
$$
\matrix{
&\ell_1^{i,J_1} =  \sum_{\lb_1} \ c_1^{i,J_1,\lb_1} \, b_{\lb_1} 
\hfill \cr\cr
&r_1^{i,J_1} \, \ell_2^{i,J_2} =  \sum_{\lb_2} \ 
c_2^{i,J_1,J_2,\lb_2} \, b_{\lb_2} \hfill \cr\cr
&\ldots \hfill \cr\cr
&r_{p-1}^{i,J_{p-1}} \, \ell_p^{i,J_p} =  \sum_{\lb_p} \ 
c_p^{i,J_{p-1},J_p,\lb_p} \, b_{\lb_p} \hfill \cr\cr
&r_p^{i,J_p} =  \sum_{\lb_{p+1}} \ c_{p+1}^{i,J_p,\lb_{p+1}} \, 
b_{\lb_{p+1}}  \, . \hfill \cr
}
$$
From (\ref{nul*}) it then follows that $\fl \,  \lb_1 , \ldots , \lb_{p+1}$
\begin{equation} \label{c*}
\sum_{i}€ \ c_1^{i,J_1,\lb_1} \cdot c_2^{i,J_1,I_2,\lb_2} \cdot 
c_{p+1}^{i,J_p,\lb_{p+1}} = 0 \, .
\end{equation}
We can conclude that
$$
\matrix{
&& \sum_i \ h_1^i \ot_{\Rc}€ \ldots \ot_{\Rc}€ h_p^i \, = \hfill \cr\cr
&&\sum_{i,J}€ \a (\ell_1^{i,J_1}) B_{J_1} \ot_{\Rc}€ 
\a (r_1^{i,J_1} \ell_2^{i,J_2})  B_{J_2} \ot_{\Rc}€ \ldots 
\ot_{\Rc}€\a (r_{p-1}^{i,J_{p-1}} \ell_p^{i,J_p})  \b (r_p^{i,J_p}) B_{J_p} 
\hfill \cr\cr
&& = \sum_{J, \lb} \sum_{i}€ \ c_1^{i,J_1,\lb_1} \ldots 
c_{p+1}^{i,J_p,\lb_{p+1}} \, \a 
(b_{\lb_1}) \, B_{J_1} \ot_{\Rc}€ \ldots 
\ot_{\Rc}€\a (b_{\lb_p}) \, \b (b_{\lb_{p+1}}) \, B_{J_p} \cr \cr
&& = 0 \, , \quad \hbox{by} \ (\ref{c*}) \, . \qquad \Box \hfill \cr
}
$$

\begin{remark} \label{full}
{\rm The algebra $\Hc_{FM}$ acts naturally on the frame bundle $F\Gb$
of every smooth \'etale groupoid $\Gb$ over $M$. Moreover, if $\Gb$ is
{\it full} , in the sense that the natural map from $\Gb$ to the 
groupoid of jets of diffeomorphisms of $M$ has dense range (i.e. 
surjects on $k$-jets for every positive integer $k$), then 
$\Hc_{FM}$ is completely determined by its action on $\Gb$. Indeed,
the proof of Proposition \ref{free} remains unchanged
if in Definition 1 one replaces the algebra
 $\, C_c^{\ify} (FM  \bar{\semi} \G_{M}) \, $
 by $\, C_c^{\ify} (\Gb)$.}  
\end{remark}

\section{The Hopf algebraic structure of ${\Hc}_{FM}$}

We devote this section to the description of the intrinsic
Hopf structure carried by the algebra $ \Hc = {\Hc}_{FM}$. 
It fits the pattern of the definition of a {\it Hopf algebroid}, 
cf. \cite[Appendix 1]{Ra} and \cite{Lu}, with
$\Hc$ as the total algebra and $\Rc$ as
the base algebra. The source and target homomorphisms
$\a: \Rc \ra \Hc$, resp. $\b: \Rc \ra \Hc$, 
obey the commutation property (\ref{ab}) and confer to
$\Hc$ its $(\Rc, \Rc)$-bimodule structure.

In order to define the coproduct, we first note that the
generators of $\Hc$ satisfy product rules when acting on $\Ac$.
Indeed, one easily checks that, for any
$a^1 , a^2 \in \Ac \,$, one has
\begin{equation} \label{pg}
\matrix{
\displaystyle \a(\ell) \, (a^1 a^2) &= &\a (\ell) (a^1) \cdot a^2 \, , 
\qquad  \fl \, \ell \in \Hc \, , \hfill \cr \cr
\displaystyle \b (r) \, (a^1 a^2) &= &a^1 \cdot \b (r) (a^2) \, , 
\qquad  \fl \,  r \in \Hc \, , \hfill \cr \cr
\displaystyle Y_i^j (a^1 a^2 ) &=& Y_i^j (a^1) \, a^2 + a^1 \, Y_i^j (a^2) 
\, , \hfill  \cr \cr
\displaystyle X_k (a^1 a^2) &=& \displaystyle X_k (a^1) \, a^2 + 
a^1 \, X_k (a^2) + \d_{jk}^i (a^1) \, Y_i^j (a^2) \, , \cr \cr
\displaystyle \d_{jk}^i (a^1 a^2) &=& \d_{jk}^i (a^1) \, a^2 + 
a^1 \, \d_{jk}^i (a^2) \, . \hfill \cr
}
\end{equation}
These identities are all of the form
\begin{equation} \label{pr}
\displaystyle h (a^1 a^2) \ = \ \sum \ h_{(1)} (a^1) \, h_{(2)} (a^2) \, , 
\quad  \hbox{with} \ h_{(1)} , h_{(2)} \in \Hc \, ,
\end{equation}
where the sum in the right hand side stands for
the customary Sweedler summation convention.
By multiplicativity, the rule (\ref{pr}) extends to
all differential operators $h \in \Hc$. While the
right hand side in (\ref{pr}) does not
uniquely determine an element $\sum \ h_{(1)} \ot h_{(2)} \in \Hc \ot \Hc$,
the ambiguity disappears in the tensor product over $\Rc$,
$\sum \ h_{(1)} \ot_{\Rc} h_{(2)} \in \Hc \ot_{\Rc} \Hc$.

\begin{proposition} \label{copr} The formula
\begin{equation} \label{D}
\D \, h = \sum \ h_{(1)} \ot_{\Rc} h_{(2)} \, , \quad \fl \, 
h \in \Hc \, ,
\end{equation}
with the right hand side given by the product rule (\ref{pr}),
defines a coproduct $\D : \Hc \ra \Hc \ot_{\Rc} \Hc$, i.e.
an $(\Rc, \Rc)$-bimodule map satisfying
\begin{itemize}
	\item[(a)] $ \D (1) = 1 \ot 1 \, ;$
	\item[(b)]
	$ (\D \ot_{\Rc} Id ) \circ \D = ( Id \ot_{\Rc} \D ) \circ \D :
	\Hc \ra \Hc \ot_{\Rc} \Hc \ot_{\Rc} \Hc \, ;$
	\item[(c)] $ \D \, (h) 
 \cdot \left( \b (b) \ot 1 \ - \ 1 \ot \a (b) \right) = 0 \, ,
	\quad \fl \, b \in \Rc \, , \,  h \in \Hc \, ,$ \\
where we have used the right action of $\Hc \ot \Hc$ on
$\Hc \ot_{\Rc} \Hc $ by right multiplication;
	\item[(d)]  
	$ \D (h_1 \cdot h_2 ) = \D (h_1) \cdot \D (h_2 ) \, , 
	\quad \quad \fl \, h_1 , h_2 \in \Hc \, ,$ \\
where the right hand side makes sense because of (c).
\end{itemize} 
\end{proposition}

\noindent {\it Proof}. The fact that the coproduct (\ref{D}) is
well-defined is an immediate consequence of 
Proposition~\ref{multid}, which also ensures that $\D $
is a bimodule homomorphism satisfying ($a$).
By the same token, the coassociativity property
($b$) is tantamount to
$$ h \, \left( (ab) \, c \right) \, = \, h \, \left( a \, (bc) \right)
\, , \quad \fl \, h \in \Hc \quad \hbox{and} \, \quad \fl \, 
a, b, c \in \Ac \, ,
$$
while ($c$) is equivalent to a special case of the above, 
corresponding to $b \in \Rc$. In turn,
($c$) ensures that the preimage of 
$\D \, (\Hc) $ in  $\Hc \ot \Hc $ is contained in the normalizer
\begin{equation} \label{Nor}
\Nc \, (\Jc) \, = \, \{ \nu \in \Hc \ot \Hc \, ; \, \nu \cdot \Jc 
\sbs \Jc \} \, ,
\end{equation}
where $\Jc$ is the right ideal generated by 
$\{ \b (b) \ot 1 \ - \ 1 \ot \a (b) \, ; \, b \in \Rc \}$
in $\Hc \ot_{\Rc} \Hc$, i.e.
such that $\Hc \ot \Hc \, / \Jc \, = \, \Hc \ot_{\Rc} \Hc$ ;
in particular, 
$$\D \, (\Hc) \sbs \Nc \, (\Jc) \,  / \Jc \, ,
$$
which is an algebra. $\qquad \Box$
\medskip

At this point let us note that the action of $\Hc$ on $\Ac$ leaves
invariant the subalgebra $\Rc_{0}€ \, = \, C_{c}€^{\ify}€ \, (FM)$,
pulled back from the space of units of the groupoid. The restriction
of $\Hc$ to End $(\Rc_{0}€)$ coincides with the algebra of differential
operators on $FM$, and therefore admits a tautological extension to
an action on $\Rc$. With this clarified, we can proceed
to define the counit.

\begin{proposition} \label{coun}
    The map $\ve: \Hc \ra \Rc$ defined by
\begin{equation} \label{eqcoun}
    \ve (h) \, = \, h(1) \, , \qquad \fl \, h \in \Hc \, ,
\end{equation}
is a bimodule map such that
 \begin{itemize}
	\item[(a)] $\ve (1) = 1 \, ;$ 
	\item[(b)] ${\rm Ker} \ve $ is a left ideal of $\Hc \, ;$
	\item[(c)]
 	$(\ve \ot_{\Rc} Id) \circ \D = ( Id \ot_{\Rc} \ve ) \circ \D = 
 	Id_{\Hc}€ \, , $ \\
where we have made the identifications
$ \Rc \ot_{\Rc}€ \Hc \, \simeq \, \Hc \, ,
\Hc \ot_{\Rc}€ \Rc \, \simeq \, \Hc \,$, the first via
the left action by $\a$ and the second via the left
action by $\b$.
\end{itemize}
\end{proposition}

\noindent {\it Proof}. The fact that $\ve: \Hc \ra \Rc$ is a bimodule
map satisfying ($a$) and ($b$) is obvious from the definition. 
In view of Proposition~\ref{multid},
($c$) amounts to
$$ 
   h \, (1 \cdot a) \, = \,  h \, (a \cdot 1) \, = \, h \, (a)
   \, , \quad \fl \, h \in \Hc \quad \hbox{and} \, \quad \fl \, 
a \in \Ac \, . 
   \hfill \qquad \Box
$$
\smallskip

The last ingredient we need is a twisted version of the antipode.
To define it, we recall that $\Ac$ carries a canonical (up
to a scaling factor) {\it faithful trace}
$\tau : \Ac \ra \Cb \, $, defined as follows
\begin{equation} \label{tr} \tau \, (f \, U_{\vp}^* )  \, = \, 
\left\{ \matrix{
\displaystyle  \int_{FM} \, f \, vol_{FM} \, , \quad \hbox{if}
\quad \vp = Id \, , \cr\cr
\displaystyle \quad 0 \, , \qquad \hbox{otherwise} \, , \hfill \cr
}\right.
\end{equation}
where $vol_{FM}$ denotes a fixed $\G_{M}$-invariant volume form
on the frame bundle, e.g. determined by the choice of a torsion-free. 
connection

\begin{proposition} \label{tap} The identity
\begin{equation} \label{ip}
\tau \, (h(a) \cdot b) = \tau \, (a \cdot \wt S (h) \, (b)) \, , \quad 
\fl \, a,b \in \Ac \, , \quad \fl \, h \in \Hc \, ,
\end{equation}
is satisfied by all operators $h \in \Hc$; it
uniquely defines an algebra anti-isomorphism ${\wt S} : \Hc \ra \Hc$,
such that
\begin{eqnarray} 
\label{ap1}	&&{\wt S}^2 \, = \, Id_{\Hc} \, , \hfill \\
\label{ap2}	&&{\wt S} \circ \b = \a \, , \hfill \\
\label{ap3}	&& m_{\Hc} \circ ({\wt S} \ot_{\Rc} Id) \circ \D = 
\b \circ \ve \circ {\wt S} : \Hc \ra \Hc \ , \hfill
\end{eqnarray}
where $m_{\Hc} : \Hc \ot \Hc \ra \Hc$ denotes the multiplication and
the left hand side of last identity makes sense because
of the preceding identity.	
\end{proposition}

\noindent {\it Proof}. Obviously, for any
$a^1 , a^2 \in \Ac \,$, one has 
\begin{equation} \label{bpg}
\matrix{
\displaystyle \tau ( \a (b) \, (a^1 ) \cdot a^2 ) &=&
 \tau ( a^1 \cdot  \b (b) \, (a^2 ) ) \, , \quad \fl \, b \in \Rc \, ,
 \hfill \cr \cr
\displaystyle \tau ( \b (b) \, (a^1 ) \cdot a^2 ) &= &
 \tau ( a^1 \cdot  \a (b) \, (a^2 ) )  \, , \quad \fl \, b \in \Rc \, .
 \hfill  \cr
}
\end{equation}
Since the fundamental vector fields are $\G_{M}$-invariant,
(\ref{pg}) and the Stokes formula gives
\begin{equation} \label{Yip}
\displaystyle \tau \, (Y_i^j (a^1) \cdot a^2 ) \, = \,
- \, \tau \, (a^1 \cdot Y_i^j (a^2)) \,
 + \, \d_{i}^{j} \, \tau \, (a^1 \cdot a^2 ) \, . \hfill
\end{equation}
On the other hand, 
for the basic horizontal vector fields, again using (\ref{pg}) 
and the Stokes formula, one gets
\begin{equation} \label{Xip}
\matrix{
 \displaystyle \tau \, (X_k (a^1) \cdot a^2 ) \,  &=& \, - \,
\tau \, (a^1 \, X_k (a^2) ) \, - \, 
\tau \, ( \d_{jk}^i (a^1) \cdot Y_i^j (a^2) ) \, \hfill \cr \cr
\displaystyle &=& \, - \, \tau \, (a^1 \, X_k (a^2) ) \, + \, 
\tau \, (a^1 \cdot \d_{jk}^i (Y_i^j (a^2) ) \, . \hfill  \cr
}
\end{equation}
Thus, the generators of $\Hc$ satisfy (\ref{ip}). 
By multiplicativity, the ``integration by part''
identity (\ref{ip}) extends to
all transverse differential operators $h \in \Hc$ and, since the trace $\tau$
is faithful, it uniquely determines the algebra
anti-homomorphism ${\wt S} : \Hc \ra \Hc$
satisfying (\ref{ap1}) and (\ref{ap2}). Finally,
(\ref{ap3}) follows from the fact that, for any $h \in \Hc$ and
$a^1 , a^2 \in \Ac \,$,
\begin{equation} 
\matrix{
\displaystyle \tau \, (a^1 \cdot a^2 \, {\wt S} (h)(1)\, )  &=& 
\tau \, (h \, (a^1 a^2 )) =
\tau \, (h_{(1)} (a^{1}) \cdot h_{(2)} (a^{2}) ) \, = \,
 \cr \cr
 \displaystyle  &=&  \tau \, (a^1 \cdot {\wt S} (h_{(1)}) h_{(2)} \, (a^2 )) 
 \, . \hfill \Box \cr
 }
\end{equation}

\begin{remark} {\rm  ${\wt S}  : \Hc \ra \Hc$ is a twisted version
of the antipode required by Lu's definition~\cite{Lu}. 
It already occurs in the flat case \cite{CM2}, 
for the genuine Hopf algebra $\Hc_n$, and will
play a similar role in the definition of the
associated cyclic module. In turn, the bimodule homomorphism
$\, \d = \ve \circ {\wt S} : \Hc \ra \Rc \, $ 
is the analogue of the modular character 
in \cite{CM2}, \cite{CM3}. }
\end{remark}

\section{Differentiable and Hopf cyclic cohomology}

The general framework for cyclic cohomology is that of 
the category of $\Lb$-modules
over the cyclic category $\Lambda$ (cf.~\cite{C3}). We recall that
the cyclic category $\Lambda$ is the small category obtained by
adjoining cyclic morphisms to the simplicial
category $\D$. The latter 
has one object $[q] = \{0 < 1 < \ldots < q \}$ for each
integer $q \geq 0$, and is generated by faces $\delta_i: [q-1] \ra 
[q]$, with $\d_{i}€$ =
the injection that misses $i$, and degeneracies $\s_j: [q+1] \ra 
[q] $, with $\s_{j}€$ =
the surjection which identifies $j$ with $j+1$, satisfying the 
relations:
\begin{equation}\label{ad}
\delta_j  \delta_i = \delta_i  \delta_{j-1}
\ \hbox{for} \ i < j  \, , \ \s_j  \s_i = 
\s_i  \s_{j+1} \qquad i \leq j 
\end{equation}
$$
\s_j  \delta_i = \left\{ \matrix{
\delta_i  \s_{j-1} \, , \hfill &i < j \hfill \cr
1_q \, , \hfill &\hbox{if} \ i=j \ \hbox{or} \ i = j+1 \cr
\delta_{i-1}  \s_j \, , \hfill &i > j+1  . \hfill \cr
} \right.
$$
The category $\Lambda$ is obtained by adjoining
for each $q$ an extra morphism $\tau_q: [q] \ra [q]$ such that
\begin{equation}\label{ae}
\matrix{
\tau_q  \delta_i = \delta_{i-1}  
\tau_{q-1} \, , &1 \leq i \leq q ,  \hfill \cr
\cr
\tau_q  \s_i = \s_{i-1} 
\tau_{q+1} \, , &1 \leq i \leq q ,\cr
\cr
\tau_q^{q+1} = 1_q \, . \hfill \cr
} 
\end{equation}

The cyclic cohomology groups of a $\Lb$-module 
$\Mc = \{ M[q] \}_{q \geq 0}€$ in the category
of vector spaces (over $\Cb$) are, by definition,
the derived groups
$$
HC^{q}€(\Mc) = Ext_{\Lambda}^{q} (\Cb , \Mc) 
\, .
$$
Using the canonical projective biresolution
of the trivial cyclic object~\cite{C3}, 
these groups can be computed as the cohomology of a bicomplex 
$ C C^{*, *} (\Mc)$ 
defined as follows
\begin{equation} \label{CbB}
\matrix{ 
&C C^{p, r} (\Mc) = \Mc [r-p], \quad r \geq p \, , \cr \cr
& C C^{p, r} (\Mc) = 0 \, ,  \quad r < p \, , \hfill \cr
}
\end{equation}
with vertical boundary operators
\begin{equation} \label{Hb}
   b = \sum_{i=0}^{q} (-1)^i \delta_i : \Mc [q-1] \ra \Mc [q] \, 
   \quad q \geq 0 \,
\end{equation}
and horizontal boundary operators
\begin{equation}\label{HB}
B =  N_q \circ \s_{-1} \circ (1_{q+1}€ - \lb_{q+1} ) : \, \Mc [q+1] \ra 
\Mc[q]  \, ,
\end{equation}
where
\begin{equation}\label{N}
 \matrix{    
    \lb_{q}€ &=& (-1)^q \tau_{q} \, , \quad
\s_{-1} = \tau_{q+1}€ \circ \s_{q}€: \Mc [q+1] \ra \Mc [q] \quad
\hbox{and} \cr \cr
      N_q &=& 1_{q}€ +  \lb_q  + \ldots +  {\lb_q}^q \, . \hfill
      }
\end{equation}
Then $H C ^{*}€(\Mc)$ is the cohomology of the first
quadrant total complex
\begin{equation} \label{TC}
    T C^{q}€(\Mc) = \sum_{p=0}^{q}€ \, C C^{p, q-p} (\Mc) \, , \label{Tot}
\end{equation}
while the cohomology of the full direct sum 
total complex 
\begin{equation} \label{PC}
    P C^{q}€(\Mc) = \sum_{p \in \Zb}€ \, C C^{p, q-p} (\Mc) ,
    \label{ST}
\end{equation}
gives the $\Zb /2$--graded {\it periodic} cyclic cohomology groups
$H C^{*}€_{\rm per}€ \, (\Mc)$.
\medskip

In particular,
the cyclic cohomology $H C^{*} (\Ac)$ 
of an associative (unital) algebra $\Ac$, 
which historically preceded the above definition (cf.~\cite{C1}), 
corresponds to the $\Lb$-module $ \Ac^{\natural}€$, with
$ \Ac^{\natural}€[q] \equiv C^{q} (\Ac)$ denoting the linear space
of ($q+1$)-linear forms $\phi$ on $\Ac$ and with $\Lb$-operators defined
as follows:
the face operators $\delta_i: C^{q-1} (\Ac) \ra C^q (\Ac) \, ,
\quad  0 \leq i \leq q \, ,$ are  
\begin{eqnarray} \label{ce}
&&\delta_i \,  \phi (a^{0}€, \ldots , a^{q}€) = 
\phi (a^{0}€, \ldots , a^{i}€a^{i+1}€ , \ldots , a^{q}€) \, ,
\quad 0 \leq i \leq q-1 \, ,\nonumber\\
&&\delta_q \, \phi (a^{0}€, \ldots , a^{q}€) = 
\phi (a^{q}€ a^{0}€,\,  a^{1}€ , \ldots , a^{q-1}€) \, ; 
\end{eqnarray}
the degeneracy operators $\s_i : C^{q+1} (\Ac) \ra C^q (\Ac) \, ,
\, \quad  0 \leq i \leq q \, ,$ are 
\begin{equation} \label{cf}
\s_i \, \phi (a^{0}€, \ldots , a^{q}€) =
\phi (a^{0}€, \ldots , a^{i}€, 1, a^{i+1}€ , \ldots , a^{q}€) 
\end{equation}
and for each $q \geq 0$ the cyclic operator   
$\tau_q : C^q (\Ac) \ra C^q (\Ac) \, $  is given by
\begin{equation}\label{cg}
 \tau_q \, \phi (a^{0}€, \ldots , a^{q}€)  =  
 \phi (a^{q}€, \, a^{0}€, \ldots , a^{q-1}€)  \, .
\end{equation}

We now take $\Ac \, = \, \Ac_{FM}$ and define its differentiable
cyclic cohomology $H C_{\rm d}€^{*}€(\Ac)$ as follows.

\begin{definition} \label{DC}
A cochain $\phi \in C^{q}€(\Ac)$ is called
 {\em differentiable} if it is of the form
 \begin{equation} \label{dc}
  \phi (a^{0}€, \ldots , a^{q}€) \, = \, \tau \, 
  \left( H( a^0 ,..., a^q ) \right)
  \, , \quad a^0 ,..., a^p \, \in {\Ac}_{FM} \, ,
 \end{equation}
where 
$$	H( a^0 ,..., a^q ) = \sum_{i=1}^r h_i^0 (a^0 ) \cdots h_i^q ( a^q )
	\, , \qquad h_{i}€^{j}€ \, \in {\Hc}_{FM} \, ,
$$
is a $q+1$-differential operator on $FM  \bar{\semi} \G_{M}$. The space of
such cochains will be denoted $ C_{\rm d}€^{q}€(\Ac) $.
\end{definition}
 
\begin{proposition} \label{Pdc}
The subspace of all differentiable cochains 
   $$ \Ac_{\rm d}€^{\natural}€ = \{ C_{\rm d}€^{q}€(\Ac) \}_{q \geq 0}
 $$
 forms a $\Lb$-submodule of $ \Ac^{\natural}€$. Furthermore, 
 for any $q \geq 1$, the map
\begin{equation} \label{Cd}
 \matrix{ 
\displaystyle \c : C^q (\Hc) \equiv
\underbrace{\Hc \ot_{\Rc} \ldots \ot_{\Rc}€ 
\Hc}_{q\hbox{-times}} \longra C_{\rm d}^{q}€(\Ac)  \cr \cr
\displaystyle  \c (h^{1} \ot_{\Rc} \ldots \ot_{\Rc}€ h^{q}€)  
  =  \tau \, (a^{0} \,  h^1 (a^1 ) \cdots h^q ( a^q )) \, ,
    \quad a^0 ,..., a^q \, \in {\Ac} \, ,
}
\end{equation}
is an isomorphism of vector spaces. 
\end{proposition} 
     
\noindent {\it Proof.} It is easy to check that the $\Lb$-operators
(\ref{ce})--(\ref{cg}) preserve the property of a cochain of being
differentiable, which proves the first claim. 
Thanks to the integration by parts property
(\ref{ip}), any differentiable $q$-cochain can be
{\it normalized}, i.e. put in the form
\begin{equation} \label{ndc}
  \phi (a^{0}€, \ldots , a^{q}€) \, = \, 
    \sum_{i=1}^r \, \tau \, (a^{0} \,  h_i^1 (a^1 ) \cdots h_i^q ( a^q ))
  \, , \quad h_{i}€^{j}€ \, \in {\Hc}_{FM} \, .
 \end{equation} 
The $q$-differential operator 
$$ 
H (a^{1}€, \ldots , a^{q}€) \, =  \,
\sum_{i=1}^r \, h_i^1 (a^1 ) \cdots h_i^q ( a^q )
$$
is uniquely determined, because of the faithfulness of the canonical
trace. Thus, the second assertion follows from Proposition \ref{multid}.
$\qquad \Box$ 

\begin{remark} \label{Rdc}
{\rm By transport of structure,
 $\Hc^{\natural}€ = \{ C^{q}€(\Hc) \}_{q \geq 0} \,$ , with 
 \begin{equation} \label{EH}
   C^{0}€(\Hc) = \Rc  \quad \hbox{and} \quad C^q (\Hc) =
\underbrace{\Hc \ot_{\Rc} \ldots \ot_{\Rc}€ 
\Hc}_{q\hbox{-times}}\, , \quad q \geq 1 \, , 
\end{equation}
acquires the structure of a $\Lb$-module.
The expressions of the
$\Lb$-operators of $\Hc^{\natural}€$ 
are virtually identical to those of the cyclic $\Lb$-module 
associated to a Hopf algebra in \cite{CM2}, with the obvious
modifications required by the replacement of the base ring $\Cb$
with $\Rc$.} 
\end{remark}
Thus, the face operators 
$\delta_i: C_{\rm d}€^{q-1} (\Ac) \ra C_{\rm d}€^q (\Ac) \, ,
\quad  0 \leq i \leq q \, ,$
are given by
\begin{equation} \label{Hce}
    \matrix{
\delta_0 (h^1 \ot_{\Rc} \ldots \ot_{\Rc} h^{q-1}) &=& 1 \ot_{\Rc} h^1 
\ot_{\Rc} \ldots \ot_{\Rc} h^{q-1} \, , \hfill \cr \cr
\delta_j (h^1 \ot_{\Rc} \ldots \ot_{\Rc} h^{q-1}) &=& 
h^1 \ot_{\Rc} \ldots \ot_{\Rc} \D h^j \ot_{\Rc} 
\ldots \ot_{\Rc} h^{q-1} \, , \cr
 &\fl& 1 \leq j \leq q-1 \, , \hfill \cr \cr
 \delta_q (h^1 \ot_{\Rc} \ldots \ot_{\Rc} h^{q-1}) &=& 
h^1 \ot_{\Rc} \ldots \ot_{\Rc} h^{q-1}
\ot_{\Rc} 1 \, , \hfill
}
\end{equation}
in particular, for $q=1$
\begin{equation}
\delta_0 (b) = 1 \ot_{\Rc}€ b = \b (b) \, , 
\quad \delta_1 (b) = b \ot_{\Rc}€ 1 = \a (b) \, 
\quad \fl \, b \in \Rc \, ;
\end{equation}
the degeneracy operators $\s_i : C_{\rm d}^{q+1} (\Ac) \ra
C_{\rm d}^q (\Ac) \, ,
\, \quad  0 \leq i \leq n \, ,$ have the expression
\begin{eqnarray} \label{Hcf}
&&\s_i (h^1 \ot_{\Rc} \ldots \ot_{\Rc} h^{q+1}) \, = \hfill \hfill 
\hfill \hfill\\
&& \qquad  h^1 \ot_{\Rc} \ldots \ot_{\Rc} h^{i} \ot_{\Rc} \ve (h^{i+1})  
 \ot_{\Rc} h^{i+2} \ot_{\Rc} \ldots \ot_{\Rc} h^{q+1} \nonumber
\end{eqnarray}
if  $\,  1 \leq j \leq q \,$ and for $q=0$
\begin{equation}
\s_0 (h) = \ve (h), \, \qquad h \in \Hc \, ; 
\end{equation}
finally, the cyclic operator 
$\tau_q : C_{\rm d}^q (\Ac) \ra C_{\rm d}^q (\Ac) \, $ is
\begin{equation}\label{Hcg}
 \tau_q (h^1 \ot_{\Rc} \ldots \ot_{\Rc} h^q)  = 
 (\D^{q-1}  \wt S (h^1)) \cdot h^2 
\ot \ldots \ot h^q \ot 1  \, ,
\end{equation}
where in the right hand side $\Hc \ot \ldots \ot \Hc$ acts on 
$\Hc \ot_{\Rc} \ldots \ot_{\Rc}€ \Hc$ by right multiplication.

The cyclic cohomology of the $\Lb$-module $\Hc^{\natural}€$ will
be denoted $H C^{*}€ (\Hc)$. By its very definition, one has
a tautological isomorphism
\begin{equation} \label{Edi}
 T^{*}€: H C^{*}€ (\Hc) \build \longra_{}^{\simeq}
 H C_{\rm d}€^{*}€(\Ac) \, .
 \end{equation}
There is a simple way for defining
the relative version of this cohomology with respect
to any compact subgroup of $K \, \sbs \,GL(n,\Rb) \, $, which
goes as follows.

For $g \in GL(n,\Rb) \,$, let $R(g)$ denote its right action on 
$\Rc = C^{\ify} (FM)$. Since the action of $GL(n,\Rb)$ on $FM$
commutes with that of
$\G_M \,$, the representation $R$ of $GL(n,\Rb)$ on $\Rc $ 
by right translations
extends to a natural action by algebra automorphisms of 
$ GL(n,\Rb) \,$ on $\, \Ac = \, C_c^{\ify} (FM  \bar{\semi} \G_{M}€ ) \, $.
Given a compact subgroup $K$, we denote by $\Ac^K$ the subalgebra of 
$K$-invariant elements in $\Ac$ and by $\, \i_K : \Ac^K \ra \Ac \, $
the corresponding inclusion map. 
We note that there is an obvious identification
$$ \Ac^K \, \simeq \,  C_c^{\ify} (FM/K  \bar{\semi} \G_{M}€ ) \, .
$$
By definition, the {\it differentiable cochains relative to} $K$
are those obtained by restricting to $\Ac^K $ the differentiable
cochains on $\Ac$,
$$  C_{\rm d}^{q}€ (\Ac , K) \, = \, 
\i_{K}^* \left( C_{\rm d}^{q}€ (\Ac) \right) 
\, .
$$
In view of Remark \ref{full}, exactly as in the case of 
$\, H C^{*}€ (\Hc) \,$, 
the complex of relative differentiable cochains
remains unchanged if $\G_M$ is
replaced by any full subpseudogroup. Because of this, we shall
denote the cyclic cohomology of the corresponding
$\Lb$-module  
\begin{equation} \label{rel}
\Ac^{K}_{\rm d}€^{\natural}€ = \{ C_{\rm d}€^{q}€(\Ac  , K ) \}_{q \geq 0}
\end{equation}
by $\, H C^{*}€ (\Hc , K) \,$ and shall refer to it as {\it 
the cyclic cohomology of $\Hc$ relative to $K$. }
\medskip 

As a simple example of an extended Hopf algebra which, without being
of the form $\Hc_{FM}€$, gives rise to a cyclic module in
a similar fashion, we shall consider the ``coarse'' extended Hopf 
algebra over 
an arbitrary associative unital algebra $\Kc$ (cf.~\cite{Lu}) 
$$ \Tc \equiv \Tc (\Kc ) := \Kc \ot \Kc^{\rm op}€ \, .
$$
The source and target maps $\, \a , \b : \Kc \ra \Tc \,$  are
$$\a (k) = k \ot 1 \, , \qquad \b (k) = 1 \ot k \, , \qquad k \in \Kc
  \, ,
$$ 
 the coproduct
  $\D : \Tc \ra \Tc \ot_{\Kc}€ \Tc \simeq \Kc \ot \Kc \ot \Kc^{\rm op}€ $
is given by
$$ \D (\ell \ot r) = (\ell \ot 1) \ot_{\Kc}€ (1 \ot r)  
  \simeq \ell \ot 1 \ot r \, , \quad \fl \, \ell, r \in \Kc  \, 
$$
the counit $\ve : \Tc \ra \Kc$  is 
$$
\ve (\ell \ot r) = \ell \cdot r \, , \quad  \ell , r \in \Kc \,
$$
and the antipode $ {\wt S} \equiv S: \Tc \ra \Tc $ is
$$ S (\ell \ot r) =  r \ot \ell \, ,
  \quad \fl \, \ell, r \in \Kc  \, . 
$$
In this case 
$$
\matrix{
C^{0}€(\Tc) &= & \Kc \, , \quad \hbox{and for} \quad
q \geq 1 \, , \hfill \cr
C^{q}€(\Tc) &= &\underbrace{\Tc \ot_{\Kc} \ldots \ot_{\Kc} 
\Tc}_{q\hbox{-times}}
    \simeq \underbrace{\Kc \ot \ldots \ot \Kc}_{q\hbox{-times}}
\ot \Kc^{\rm op} \, . \cr
}
$$
With the latter identification, the face operators
$\delta_i: C^{q-1} (\Tc) \ra C^q (\Tc) \, ,$
$  0 \leq i \leq q \, ,$ are given by
$$
 \matrix{
\d_0 (k^1 \ot \ldots \ot k^q) = 1 \ot k^1 \ot \ldots \ot k^q \, ,
\quad  k^1 , \ldots , k^q  \in \Kc \, ,\hfill \cr\cr
\d_i (k^1 \ot \ldots \ot k^q) = k^1 \ot \ldots \ot k^{i} \ot 1 \ot k^{i+1} 
\ot \ldots \ot k^q \, ,
\quad i = 1 , \ldots , q-1 \hfill \cr\cr
\d_q (k^1 \ot \ldots \ot k^q) = k^1 \ot \ldots \ot k^q \ot 1 \hfill \cr
} 
$$
for $q \geq 1$ while for $q=0$,
$$
\d_0 (k) = 1 \ot k \ , \qquad \d_1 (k) = k \ot 1 \, ;
$$
the degeneracy operators 
$\s_i : C^{q+1} (\Tc) \ra C^q (\Tc) \, ,
\, \quad  0 \leq i \leq q \, , \quad$ become
$$
\s_i (k^1 \ot \ldots \ot k^{q+2}) = k^1 \ot \ldots \ot k^{i+1} \, k^{i+2}
\ot \ldots \ot k^{q+2} \, ;
$$
finally, the cyclic operator $\tau_q : C^q (\Tc) \ra C^q (\Tc) \, $ is 
exactly the cyclic permutator
$$
\tau_{q}€ (k^1 \ot \ldots \ot k^{q+1}) = 
k^2 \ot \ldots \ot k^{q+1} \ot k^{1} \, .
$$

\begin{lemma}  \label{triv}
    For any associative unital algebra $\Kc$, one has
$$ H C^{*}€(\Tc (\Kc)) \, \simeq  H C^{*}€(\Tc (\Cb)) \, 
\equiv H C^{*} (\Cb) \, , 
$$ 
the isomorphism
being induced by the unit map $\eta : \Cb = \Tc (\Cb) \ra \Tc (\Kc)$.
\end{lemma} 

\noindent {\it Proof.} Fix a linear functional $\nu \in \Ac^* \,$ with
$\nu (1) = 1 \,$ and define 
the operator $s : C^{q}€(\Tc) \ra C^{q-1}€(\Tc)$ by setting
$$
s (k^1 \ot \ldots \ot k^{q+1}) = \nu (k^1) \, k^2 \ot \ldots \ot k^{q+1}
\, , \quad k^1 , k^2 , \ldots , k^{q+1} \in \Kc \, .
$$
It can be easily checked that
$$
s \, \d_i = \d_{i-1} \, s \, , \qquad i = 1,\ldots , q \, .
$$
Thus, the $\Lb$-module map $\eta^{\natural}€ : \Cb^{\natural}€
\ra \Tc^{\natural}€$, associated to the unit map, induces an 
isomorphism for Hochschild and  therefore for cyclic cohomology
too. $\qquad \Box$
\medskip

This lemma will be used to show that, in the case of a flat
$n$-manifold, the cyclic
cohomology of the extended Hopf algebra coincides with that of
the Hopf algebra $\Hc_n$ introduced in \cite{CM2}. We preface
the statement with the
remark that $\Hc_n$ can be defined as the subalgebra of
$\Hc_{F\Rb^{n}€}€$ generated over $\Cb$ by the operators
$\{X_{k}€, Y_{i}€^{j}€, \d_{jk}€^{i}€\}$ corresponding to the
trivial flat connection of $\Rb^{n}€$. 

\begin{proposition}
If $N^n$ is a flat affine manifold, 
then the canonical inclusion $\kappa : \Hc_n \ra \Hc_{FN} $ 
associated to the flat connection, induces isomorphism in cyclic
cohomology
$$
\kappa^K_* : HC^* (\Hc_n , K) \build \longra_{}^{\simeq} HC^* (\Hc_{FN} , K) \, ,
$$
for any compact
subgroup $K \, \sbs \,GL(n,\Rb) $
\end{proposition}

\noindent {\it Proof.} As before, we shall use the abbreviated
notation
$$ \Hc = \Hc_{FN} \, , \quad \Ac = \Ac_{FN}€ \quad \hbox{and} \quad
  \Rc = \Rc_{FN}€ \, . 
$$
Also, we shall identify $\Hc_n$ with its image via the 
homomorphism of extended Hopf algebras
$\kappa : \Hc_n \ra \Hc $ 
induced by the Hopf action of $\Hc_n$ on $\Ac$.

Applying Proposition \ref{free}, one gets a canonical
isomorphism of $(\Rc, \Rc)$-bimodules
\begin{equation} \label{cog}
\Hc \simeq \a (\Rc) \ot \b (\Rc) \ot \Hc_n \, 
\simeq \Tc \ot \Hc_n \, ,
\quad \hbox{where} \quad \Tc = \Rc \ot \Rc^{\rm op}€ \, .
\end{equation}
Furthermore, one can easily check that 
$$ \D_{\Hc}€ = \D_{\Tc}€ \ot \D_{\Hc_n}€ \quad \hbox{and} \quad
  \ve_{\Hc}€ = \ve_{\Tc}€ \ot \ve_{\Hc_n}€\, ,
$$
that is, as coalgebroid, $\Hc$
is actually isomorphic to the external tensor product between
the coarse coalgebroid $\Tc$ over $\Rc$ and the coalgebra $\Hc_n$
over $\Cb$.
This implies that 
$$
  \d_{i}€^{\Hc}  = \d_{i}€^{\Tc} \ot \d_{i}€^{\Hc_n} \, \quad \hbox{and} 
  \quad
        \s_{j}€^{\Hc} = \s_{j}€^{\Tc} \ot \s_{j}€^{\Hc_n}  \, .
$$
Thus, in the category of cosimplicial modules, one has
$$
\Hc^{\natural} \simeq \Tc^{\natural} \ts \Hc_n^{\natural} \, .
$$
Applying the Eilenberg-Zilber theorem it follows that
\begin{equation} \label{EZ}
    HH^* (\Hc) \simeq HH^{*}€(\Tc) \ot HH^* (\Hc_n) \simeq HH^* 
    (\Hc_n) \, ,
\end{equation}
with the second isomorphism being a consequence of Lemma \ref{triv}. 
By the functoriality of Eilenberg-Zilber isomorphism, 
the composition of the two
isomorphisms in (\ref{EZ}) is induced by the canonical homomorphism
$\kappa : \Hc_n \ra \Hc $, that is
\begin{equation} \label{HHi}
\kappa_* : HH^* (\Hc_n) \build \longra_{}^{\simeq} HH^* (\Hc) \, .
\end{equation}
But $\kappa : \Hc_n \ra \Hc $ is a homomorphism of extended Hopf 
algebras and therefore gives rise to
a morphism of cyclic modules 
$\kappa^{\natural}€ : \Hc_n^{\natural}€ \ra \Hc^{\natural}€ $.

The stated isomorphism for the absolute case now follows from
(\ref{HHi}) and the exact sequence relating Hochschild and
cyclic cohomology. The relative case is proved by restriction to
$K$-invariants. $\qquad \Box$    
\medskip

We are now in a position to prove the first main result
of this paper, asserting that $H C^* (\Hc_{FM})$ and its
relative variants depend only
on the dimension $n$ of the manifold. In view of with 
\cite[\S 7, Theorem 11]{CM2},
which identifies $H C^* (\Hc_{n})$ to the Gelfand-Fuchs cohomology, this
result provides a cyclic analogue 
to Haefliger's Theorem IV.4 in \cite{H2}.

\begin{theorem} \label{MI}
For any $n$-dimensional manifold $M$ and 
for any compact subgroup $K \, \sbs \,GL(n,\Rb) $, 
one has a canonical isomorphism
$$   H C^* (\Hc_{FM} , K) \, \simeq \, H C^* (\Hc_{n} , K ) \, .
$$
\end{theorem} 

\noindent {\it Proof.} To construct the stated isomorphism, we shall
have to specify certain Morita equivalence data. The resulting
isomorphism however will be independent of the choices made. 

Concretely, we fix an open cover of $M$ 
by domains of local charts $\Uc = \{ V_i \}_{1 \leq i \leq r} $
together with a partition of unity 
subordinate to the cover $\Uc \, ,$
$$
\sum \, \chi_{i}^{2}€ (x) = 1 \, , \quad \chi_{i} \in C_c^{\ify} (V_{i}) \, .
$$
We begin by forming the smooth \'etale groupoid
$$
  G_{\Uc} = \{ (x, i, j) ; \ x \in V_i \cap V_j \, , \quad
  1 \leq i, j \leq r \} \, ,
$$
whose space of units is the flat manifold 
$N \, = \, \coprod_{i=1}^r \, V_{i} \ts \{ i \} \, .$
One has a natural Morita equivalence between the algebras
$\, C_c^{\ify} (M) \, $ and $\, C_c^{\ify} (G_{\Uc}) \,  $,
and a similar construction
continues to function in the presence of the 
pseudogroup $\G_M$, as well as at the level of the frame bundle.

More precisely, there is a {\it full pseudogroup} $\G_{\Uc}$ on
$N$ such that the corresponding smooth \'etale groupoid 
$N \bar{\semi} \G_{\Uc}$ contains $G_{\Uc}$ and such that the algebra
$ \Ac =  C_c^{\ify} (FM \bar{\semi} \G_M ) \,$ can be identified with
the reduction by a canonical idempotent of the groupoid algebra
$ \Bc = C_c^{\ify} (FN \bar{\semi} \G_{\Uc}) \,$,
\begin{equation} \label{Ide}
  \Ac \, \sm \,
	e \, \Bc \, e  \, , \quad e^{2}€ = e \in \Bc \, ;
\end{equation}
at the same time, one also has a canonical identification
$$ C^{\ify} (FM) \,  \sm \, e \, C^{\ify} (FN) \, e \, .
$$

Indeed, the elements of 
$\, \Bc \, = \, C_c^{\ify} (FN \bar{\semi} \G_{\Uc}) \,$ 
can be represented as finite sums of the form
\begin{equation} \label{B}
b \, = \, \sum_{\vp \in \G_{M}€}€ \, \sum_{i, j}€ \,
f_{\vp,ij} \, U_{\vp_{ij}}^{*}  \, , \quad
f_{\vp,ij} \in C_c^{\ify} ({\rm Dom} \vp_{ij} ) \, ,
\end{equation}
where, for any $\vp \in \G_{M} $, 
$$\vp_{ij}: 
\vp^{-1} \, \left( \vp \, ({\rm Dom} \vp \, \cap \, V_{i} ) \, \cap 
 V_{j} \right) \ts \{i\} \,  \ra \,
\vp \, ({\rm Dom} \vp \, \cap \, V_{i} ) \, \cap V_{j} \ts \{j\}
$$
stands for the obvious identification given by the
restriction of $\vp$. The definition of the 
multiplication is as follows:
\begin{equation} \label{mm}
\psi_{jk} \circ \vp_{ij} \, = \, (\psi \circ \vp)_{ik} \quad 
\hbox{or equivalently} \quad 
  U_{\vp_{ij}}^{*} \cdot U_{\psi_{jk}}^{*} =
  U_{(\psi \circ \vp)_{ik}}^{*} \, ;
\end{equation}
when $\vp = Id$, instead of $ U_{\vp_{ij}}^{*}$ 
we shall simply write $U_{ij}^{*}$. 
The canonical trace
${\wt \tau} : \Bc \ra \Cb \, $ is related to the trace
$\tau : \Ac \ra \Cb \, $ by the formula
\begin{equation} \label{Tr}
{\wt \tau} \left( \sum_{\vp \in \G_{M}€}€ \, \sum_{i, j}€ \,
f_{\vp,ij} \, U_{\vp_{ij}}^{*} \right) \, = \,
\tau \left( \sum_{\vp \in \G_{M}€}€ \, \sum_{i}€ \,
f_{\vp,ii} \, U_{\vp_{ii}}^{*} \right) \, .
\end{equation}

With the above notation, the idempotent stipulated by (\ref{Ide})
has the expression
$$
e \, = \, \sum_{i,j}€ {\wt \chi}_{i} \, U_{ij}^{*} \, {\wt 
\chi}_{j}  \, ,
\quad \hbox{with} \quad {\wt \chi}_{i} = \chi_{i} \circ \pi \, 
$$
and the claimed identification is given by the map
\begin{equation} \label{Mor}
  \i  (f \, U_{\vp}^{*}€) \, = \, \sum_{i,j}€ \, {\wt \chi}_{i}€ \, 
  f \, U_{\vp_{ij}€}^{*}€ \, {\wt \chi}_{j}€  \, , 
  \quad f \in C_c^{\ify} (FM) \, , \, \vp \in \G_{M}€ \, .
\end{equation}
The corresponding pair of bimodules which implements the Morita
equivalence between $\Ac$ and $\Bc$ consists of 
$$ \Pc  =  e \Bc \quad \hbox{and} \quad \Qc = \Bc e \, ,  \quad
\hbox{with} \quad \Pc \ot_{\Bc} \Qc \simeq \Ac \quad \hbox{and}
\quad \Qc \ot_{\Ac} \Pc \simeq \Bc \, .
$$
Setting for any $ \,  i \, = \, 1, ..., r \, $,
$$ \matrix{
u_i = \sum_k {\wt \chi}_{k} U_{ki}^* \cdot {\wt \chi}_{i} \, \in \Pc 
\, , \quad v_i = {\wt \chi}_{i} \cdot \sum_k  U_{ik}^* {\wt \chi}_{k} \,
\in \Qc \, , \quad \hbox{respectively} \hfill \cr \cr
 \quad u'_i = \sum_k {\wt \chi}_{k} U_{ki}^* 
\, , \quad e \cdot u'_i = u'_i \, \quad \hbox{and} \quad
\quad v'_i = \sum_k  U_{ik}^* {\wt \chi}_{k} \, , \quad
v'_i \cdot e = v'_i \, ,
}
$$
one has
$$ \sum_i \, u_i \, v_i \, = \, e \, , \quad \hbox{and}
\quad \sum_i \, v'_i \, u'_i \, = \, \sum_i U_{ii}^* \, .
$$
The constructive proof of the Morita invariance for the
Hochschild and cyclic cohomology of algebras (cf. \cite{Mc})
associates to these data canonical cochain equivalences
$\Ps^{*}€ = \{ \Ps^{q}€ \} \,$
and $\, \Th^{*}€ = \{ \Th^{q}€ \} \, $, which are homotopic
inverses to each other. 
We shall only have to check that 
these cochain maps, as well as the cochain homotopies, 
preserve the corresponding subcomplexes of differentiable
cochains. 

The cochain map 
$\Ps^{*}€ : C^{*}€(\Ac) \ra C^{*}€(\Bc) \,$ is given by
$$ \matrix{ 
\Ps^{q}€(\ph) \, \left( \sum_{i,j}€ \, f^{0}_{ij}€ \,
U_{\vp^{0}_{ij}€}^{*}€ \, , \sum_{i,j}€ \, f^{1}_{ij}€ \,
U_{\vp^{1}_{ij}€}^{*}€ \, , \ldots ,
\sum_{i,j}€ \, f^{q}_{ij}€ \, 
U_{\vp^{q}_{ij}€}^{*}€  \right) \, = \hfill \cr \cr
\sum_{i_{0}€, i_{1}€, \ldots , i_{q}€} \,  \ph \, \left( 
 f^{0}_{i_0 i_1} \, 
U_{\vp^{0}}^{*}€ \, , \, f^{1}_{i_1 i_2} \, U_{\vp^{1}}^{*}€  
\, , \ldots , \, f^{q}_{i_q i_0} \, U_{\vp^{q}}^{*}€ \, \right) \, . \cr
}
$$
Assuming $\ph \in C_{\rm d}^{q}€ (\Ac) \,$, of the form
\begin{equation} \label{cd}
\phi (a^{0}€, \ldots , a^{q}€) \, = \, 
   \tau \left( h^{0}€ (a^{0}) \cdot h^1 (a^1 ) \cdots h^q ( a^q )
   \right) \, , \quad \hbox{with} \quad h^{i}€ \in \Hc_{FM}€ ,
\end{equation}
one has
$$ \matrix{ 
\Ps^{q}€(\ph) \, \left( \sum_{i,j}€ \, f^{0}_{ij}€ \, 
U_{\vp^{0}_{ij}€}^{*}€ \, , \sum_{i,j}€ \, f^{1}_{ij}€ \,
U_{\vp^{1}_{ij}€}^{*}€ \, ,
 \ldots , \sum_{i,j}€ \, f^{q}_{ij}€ \, U_{\vp^{q}_{ij}€}^{*}€ 
\right) \, = \hfill \cr \cr
\sum_{i_{0}€, i_{1}€, \ldots , i_{q}€}€ \, \tau€ \left(  
h^{0}€ ( f^{0}_{i_{0}€\,i_{1}€} \, U_{\vp^{0}}^{*}€ ) \,
h^{1}€ ( f^{1}_{i_{1}€\,i_{2}€} \, U_{\vp^{1}}^{*}€ ) \,
\ldots \, h^{q}€ ( \, f^{q}_{i_{q}€\,i_{0}€} \, 
U_{\vp^{q}}^{*}€ ) \right) \, ; 
}
$$
using (\ref{Tr}) it is easily seen that 
the latter expression is of the form
$$ €
\sum_s \, {\wt \tau}€ \left( \sum_{i,j}€ {\wt h}^0_{s} 
(f^{0}_{ij} U_{\vp^{0}_{ij}€}^{*}€)
\, \cdot \, \sum_{i,j}€ {\wt h}^1_{s} (f^{1}_{ij} U_{\vp^{1}_{ij}€}^{*}€) \,
\cdots \,
\sum_{i,j} {\wt h}^q_{s} (f^{q}_{ij} \, U_{\vp^{q}_{ij}€}^{*}€ ) \right) \, , 
$$
with $\, {\wt h}^0_{s} \, , \ldots , {\wt h}^q_{s} $  transverse
differential operators on $FN \bar{\semi} \G_{\Uc}$.

The cochain map $\Th^{*}€ : C^{*}€(\Bc) \ra C^{*}€(\Ac) \,$ has the
expression
$$ \matrix{ 
\Th^{q}€({\wt \ph}) \, \left( f^{0}€ \, U_{\vp^{0}} \, , 
f^{1}€ \, U_{\vp^{1}} \, , \ldots ,
 f^{q}€ \, U_{\vp^{q}} \right) \, =  \hfill \cr \cr
 {\wt \ph} \, \left(  \sum_{i_{0}€, i_{1}€, \ldots , i_{q}€}€ \,
{\wt \chi}_{i_0}€
  f^{0}€ U_{\vp^{0}_{i_0 i_1}€}^{*}€ {\wt \chi}_{i_1}€ , 
 {\wt \chi}_{i_1}€  f^{1}€ \, U_{\vp^{1}_{i_1 i_2}€}^{*}€ {\wt \chi}_{i_2}€ , 
    \ldots , {\wt \chi}_{i_q}€
 f^{q}€ U_{\vp^{q}_{i_q i_0}€}^{*}€ {\wt \chi}_{i_0}€  \right)
\, . \cr
}
$$
When ${\wt \ph} \in C_{\rm d}^{q}€ (\Bc) \,$ is of the form
\begin{equation} \label{tcd}
{\wt \phi} (b^{0}€, \ldots , b^{q}€) \, = \, 
   {\wt \tau}€ \left( {\wt h}^{0}€ (b^{0})  {\wt h}^1 (b^1 ) \cdots 
   {\wt h}^q ( b^q )
   \right) \,  \quad \hbox{with} \quad {\wt h}^{i}€ \in \Hc_{FN}€ 
   \, ,
\end{equation}
one gets
$$ \matrix{ 
\Th^{q}€({\wt \ph}) \, \left( f^{0}€ \, U_{\vp^{0}} \, , 
f^{1}€ \, U_{\vp^{1}} \, , \ldots ,
 f^{q}€ \, U_{\vp^{q}} \right) \, =  \hfill \cr \cr
\sum_{i_{0}€, i_{1}€, \ldots , i_{q}€}€ \, {\wt \tau}€ \left(
{\wt h}^{0}€ ({\wt \chi}_{i_{0}€} f^{0}
U_{\vp^{0}_{i_{0}€\,i_{1}€}}^{*}€ {\wt \chi}_{i_{1}€} ) \cdot
{\wt h}^{1}€ ({\wt \chi}_{i_{1}€} f^{1} 
U_{\vp^{1}_{i_{1}€\,i_{2}€}}^{*}€ {\wt \chi}_{i_{2}€} )
\cdots
{\wt h}^{q}€ ({\wt \chi}_{i_{q}€} f^{q}
U_{\vp^{q}_{i_{q}€\,i_{0}€}}^{*}€ {\wt \chi}_{i_{0}€} ) 
\right) ;
}
$$
again, one can recognize
the last expression to be of the form
$$ \sum_{s}€
\tau \left( h^0_{s}€ (f^{0} \, U_{\vp^{0}}^{*}€) \cdot
h^1_{s}€ (f^{1} \, U_{\vp^{1}}^{*}€)
\cdots  h^q_{s}€ (f^{q} \, U_{\vp^{q}}^{*}€ ) \right) \, , 
$$
with $\,  h^0_{s}€ , h^1_{s}€ , \ldots , h^q_{s}€
\in \Hc_{FM}€ $.

Similar arguments apply to the canonical cochain homotopies used in
showing that $\Ps^{*}€$ and $\Th^{*}€$ are homotopic inverses to
each other in Hochschild cohomology.
Since both $\Ps^{*}€$ and $\Th^{*}€$  are actually cochain maps for 
the cyclic bicomplex, it follows
that they induce the stated isomorphism. 
All these constructions
are clearly equivariant with respect to $O(n)$ and therefore
apply, by restriction to $K$-invariants,
to the relative cohomology as well.  $\qquad \Box$
\medskip

We conclude this section with the observation that the ordinary
(de Rham) homology of the manifold $M$ can also be construed as 
differentiable (or Hopf) cyclic cohomology. The extended Hopf algebra
responsible for this interpretation is the algebra
of the differential operators on $FM$,
$$ \Dc \, := \, \Dc_{FM} \, ,
$$
viewed as the algebra of linear transformation of 
$\Rc_0 = C_c^{\ify}(FM)$
generated by the vector fields on $FM$, acting as derivations,
and by the 
functions in $\Rc = C^{\ify} (FM)$, acting as multiplication  
operators via $\a \equiv \b$. Like its extension $\Hc$, 
$\Dc$ also carries 
an intrinsic Hopf-algebraic structure in the sense of \S 2. The
coproduct $\D : \Dc \ra \Dc \ot_{\Rc} \Dc $ is given by the 
analogue of (\ref{D}), determined by the 
Leibniz rule (\ref{pr}) corresponding to its action on $\Rc$. 
Similarly, 
the twisted antipode $\wt{S} : \Dc \ra \Dc$ is determined by
the analogue of
integration by parts formula (\ref{ip}), with respect to
 $\tau \vert \Rc_0 $ which is given by the
integral against the canonical volume form on $FM$. In other words,
$\Dc$ inherits the quotient Hopf structure corresponding to
the restricted action of $\Hc$ on $\Rc$. 

The differentiable cyclic
cochains on $\Rc_0$ are defined by specializing Definition \ref{DC}
to the algebra $\Rc_0$, acted upon by $\Dc$. By
obvious analogues of Proposition \ref{Pdc}
and Remark \ref{Rdc}, the corresponding $\Lb$-module can also be described
purely in terms of the Hopf structure of $\Dc$. The ensuing cyclic
complex, unlike that of $\Dc$ viewed solely as an
algebra (cf. \cite{BG}, \cite{Wo}, for the cyclic homology of the latter),
is quasi-isomorphic to the full cyclic
complex of the algebra $\Rc_0$. In a sense made more precise by the proof
below, the differentiable cyclic cochains are in the same relation to
the continuous cyclic cochains as the smooth currents are in relation
to arbitray currents on $FM$.

\begin{proposition} \label{mi}
The tautological action
of $\Dc_{FM}$ on $C_c^{\ify}(FM)$ or, equivalently, the 
inclusion of the differentiable cyclic subcomplex 
into the usual cyclic complex of $C_c^{\ify}(FM)$, induces 
isomorphism in cyclic cohomology, resp. relative cyclic cohomology, 
$$ \matrix{
 H C_{\rm per}^* (\Dc_{FM}) \, &\simeq& \,  
H C_{\rm per}^* \left(C_c^{\ify}(FM) \right) \, , \quad 
\hbox{resp.} \hfill \cr \cr
H C_{\rm per}^* (\Dc_{FM}, O(n)) \, &\simeq& \,  
H C_{\rm per}^* \left(C_c^{\ify}(PM) \right) \, , \qquad
 PM := FM/O(n) .
}
$$
\end{proposition}

\noindent {\it Proof.} Since the map that sends 
$\vf \in \Om^{{\rm dim} FM \, - \, q} \, (FM)$ 
to the cyclic cochain
on $C_c^{\ify}(FM)$
$$ \g_{\vf} \, (f^0 , f^1 , \ldots , f^q )
\, = \, \int_{FM} \, f^0 \, d f^1 \wdg \ldots  \wdg  df^q \wdg \vf
$$
induces isomorphism between $\sum_{ i \equiv *} \, H^{n(n+1)-i} (FM)$ 
and $H C_{\rm per}^* \left(C_c^{\ify}(FM) \right)$, it suffices to
show that any such cochain $\g_{\vf}$ is differentiable. To check this, 
we fix a torsion-free connection $\om$ 
on $FM$ and let
$\{ X_k , Y_i^j \} $ (resp. $\{\t^k , \om^i_j \}$) stand for the 
corresponding basis of the tangent (resp. cotagent) space to $FM$
at any point. Then, for each $\ell = 1, \ldots , q \, ,$ one has
$$
df_{\ell} =  Y_i^j (f_{\ell}) \, \om_j^i +  X_k (f_{\ell}) \, \t^k \, . 
$$ 
On the other hand, $\vf$ can be expressed as a linear combination
over $\Rc$ of 
monomials formed from the basis $\{\t^k , \om^i_j \}$. Thus,
$$ \g_{\vf} \, (f^0 , f^1 , \ldots , f^q )
\, = \, \sum_s \, \int_{FM} \,  f^0 \, h_s^1 (f^1 ) 
 \cdots   h_s^q (f^q ) \, vol_{FM} \, , \quad 
 \hbox{with} \quad h_s^{\ell} \in \Dc.
$$
The relative case can be proved in a similar fashion. 
$\qquad \hfill \Box$

\section{Geometric realization of the cyclic van Est isomorphism}   

As mentioned before, in conjunction with \cite[\S 7, Theorem 11]{CM2}, 
Theorem \ref{MI} implies
that, for any $n$-dimensional smooth manifold $M$,
$\, H C^{*}€ (\Hc_{FM}€)$ is canonically isomorphic to
the Gelfand-Fuchs cohomology
groups $H^*({\Fa}_n)$. The purpose of this section is to give
a geometric construction of an
explicit cochain map, from the Lie algebra cohomology complex
to the cyclic bicomplex, which implements this cyclic analogue
of the van Est isomorphism.
As in the flat case \cite[\S 7]{CM2}, the cochain map
will be assembled in two stages.
The first consists in mapping the Gelfand-Fuchs cohomology
of the Lie algebra ${\Fa}_{n}$ of formal vector
fields on $\Rb^n$ to the 
$\Gc_M$-equivariant cohomology of $FM$, where $\Gc_M := {\rm Diff} (M)$. 
The second involves the
canonical map $\Phi$ relating the equivariant cohomology
to the cyclic cohomology \cite[Theorem 14, p.220]{Co}.

In order to construct the first map,
we shall fix a torsion-free linear connection $\nabla$ on $M$,
with connection form $\om = (\om_j^i)$. We denote by
$F^{\ify}M$ the bundle of frames of infinite order on $M$,
formed of jets of infinite order $j_0^{\ify} (\psi)$ 
of local diffeomorphisms $\psi$,
with source $0$, from $\Rb^n$ to $M$, and by $\pi_{1}:
F^{\ify}M \ra FM$ its projection to $1$-jets.
The connection $\nabla$ determines a cross-section 
$$\s \equiv \s_{\nabla} : FM \ra F^{\ify}M \, , \quad
\pi_1 \circ \s = {\rm Id} \, ,
$$
by the formula
\begin{equation} \label{eq4.1}
\s (u) = j_0^{\ify} (\exp_x^{\nabla} \circ u) \ , \qquad u \in F_x M \, .
\end{equation}
This cross-section is clearly ${\rm GL}(n, \Rb)$-equivariant
\begin{equation} \label{eq4.2}
\s_{\nabla} \circ R_a = R_a \circ \s_{\nabla} \, , \qquad 
 a \in {\rm GL} (n, \Rb) \, ,
\end{equation}
as well as Diff-equivariant, 
\begin{equation} \label{nat}
\s_{\nabla^\vp} = \wt{\vp}^{-1} \circ \s \circ \wt{\vp} \ , 
\qquad \fl \, \vp \in \Gc_M  \, ;
\end{equation}
here $\nabla^\vp$ corresponds to $\wt{\vp}^* \om$, 
where by $\wt \vp$ we denote the action of
$\vp$ on both $FM^{\ify}$ and $FM$.

For each simplex $(\vp_0 , \ldots , \vp_p) \ts \D^p$, we define the 
 map
$$
\s_{\nabla} (\vp_0 , \ldots , \vp_p) : \D^p \ts FM \ra F^{\ify}M
$$
by the formula
\begin{equation} \label{eq4.5}
\s_{\nabla} (\vp_0 , \ldots , \vp_p) (t,u) = 
\s_{\nabla (\vp_0 , \ldots , \vp_p ; t)} (u)
\end{equation}
where
$$
\nabla (\vp_0 , \ldots , \vp_p ; t) = \sum_0^p \ t_i \, \nabla^{\vp_i} \, .
$$
Note that $\s_{\nabla (\vp_0 , \ldots , \vp_p ; t)}$ depends 
only on the class $[\vp_0 , \ldots 
, \vp_p ; t] \in E\Gc_M$, where $E \Gc_M $ is the geometric
realization of the simplicial set $N\Gc_M$, with 
$$ N\Gc_M [p] \, = \, \underbrace{\Gc_M \times \, \cdots \, \times 
\Gc_M}_{p+1\hbox{-times}} \, .
$$
Also, in view of (\ref{nat}),
one has
\begin{equation} \label{ic}
\s_{\nabla} (\vp_0 \, \vp , \ldots , \vp_p \, \vp) (t,u) = \wt{\vp}^{-1} \, 
\s_{\nabla} (\vp_0 , \ldots , 
\vp_p) (t, \wt{\vp}(u)) \, .
\end{equation}

At this point, let us recall that the $\Gc_M$-equivariant cohomology of
$FM$, twisted by the orientation sheaf, can be computed by means
of the bicomplex 
$$ \left( C^{\ast, \ast} (\Gc_M ; FM) , \d , \partial \right) \, ,
$$
defined as follows:
$C^{p, m} (\Gc_M ;  FM) = 0$ unless 
 $p \geq 0$ and $ - {\rm dim \,} FM \leq m \leq 0$, when
 $ C^{p, m} (\Gc_M ;  FM) \equiv C^p (\Gc_M , \Om_{-m} \,( FM)) $
is the space of totally antisymmetric maps $\mu : \Gc_M^{p+1}
\ra \Om_{-m} \,( FM))$, from $\Gc_M^{p+1}$ to the currents of
dimension $-m$ on $FM$, such that
\begin{equation} \label{IC}
\mu (\vp_0 \, \vp , \ldots , \vp_p \, \vp) = ({\wt \vp}_* )^{-1} \, \mu
(\vp_0 , \ldots , \vp_p) \, , 
\end{equation}
with ${\wt \vp}_*$ denoting the transpose of ${\wt \vp}^*$ acting on
the forms $\Om (FM)$; the operator $\d$ is the simplicial
coboundary and $\partial$ is the de Rham boundary for currents.

We denote by $C^* ({\Fa}_n)$ the Lie algebra cohomology complex of the
antisymmetric multilinear functionals on the Lie algebra ${\Fa}_{n}€$
of formal vector fields on $\Rb^n$, which are continuous with respect to the
$\Ic$-adic topology,  i.e. depend only on finite jets at 
$0 \in \Rb^n$ of vector fields. 
The canonical flat connection of 
$F^{\ify}M$ determines an isomorphism between $C^* ({\Fa}_n)$ and
the space $\Om^{*}€ \, ((F^{\ify}€M))^{\G_{M}€}€$
of all $\G_M$-invariant forms on $F^{\ify}M$; we shall denote by
${\wt \vf}$ the $\G_M$-invariant form
corresponding to $\vf \in C^* ({\Fa}_n)$.

With this notation in place, we now define for any 
$\vf \in C^{q}€ ({\Fa}_n)$ and
for any pair of 
integers $(p, m) \, $ such that $ \, p \geq 0 , \, 
 - n(n+1) \leq m \leq 0 \,$
and $\, p + m = q - n(n+1)$, a current of dimension $-m$ on $FM$ by
the following formula, where $\, \eta \in \Om_{c}€^{-m}€ \, (FM) $,
\begin{equation} \label{dchain}
\matrix{
\langle C_{p,m} (\vf) (\vp_0 , \ldots , \vp_p) , \eta \rangle && =
   \hfill \cr\cr
&& \displaystyle (-1)^{\frac{m(m+1)}{2}}€
\int_{\D^p \ts FM} 
\s_{\nabla} (\vp_0, \ldots, \vp_p)^* (\pi_1^* (\eta) \wedge \wt{\vf})
 \hfill \cr\cr
&& = \displaystyle (-1)^{\frac{m(m+1)}{2}}€
\int_{\D^p \ts FM}  \eta \wedge \s_{\nabla} (\vp_0 , 
\ldots , \vp_p)^* (\wt{\vf}) \, .
\hfill \cr
}
\end{equation}

\begin{lemma} \label{Lchain}
    For any   $\vf \in C^{q} ({\Fa}_n)$, one has
    $\, C_{p,m} (\vf) \in C^p (\Gc_M , \Om_m \,( FM)) \,$ and
the assignment
$$ C(\vf) = \sum_{p+m=q-n(n+1)}€\,  C_{p,m} (\vf) \, ,
$$
defines a map of (total) complexes
\begin{equation} \label{chain}
    C: \left( C^* ({\Fa}_n) , \, d \right) \, \ra \,
\left( T C^{*} (\Gc_M ;  FM) , \, \d + \partial \right) \, .
\end{equation}
\end{lemma}

\noindent{\it Proof.} We first check the identity
\begin{equation} \label{eq4.8}
C_{p,m} (\vf) (\vp_0 \, \vp , \ldots , \vp_p \, \vp) = \vp_*^{-1} \,
C_{p,m} (\vf) (\vp_0 , \ldots , \vp_p) \, .
\end{equation}
Indeed, using (\ref{ic}), one has $\, \fl \eta \in  \Om_c^{-m} \,( FM)) \,$
$$
\matrix{
&&\displaystyle \int_{\D^p \ts FM} \s_{\nabla}€ (\vp_0 , \ldots , 
\vp_p)^* \, \left( {\wt \vp}^{-1\, *} (\pi_{1}€^{*}€(\eta)) \wedge
{\wt \vp}^{-1\, *} (\wt{\vf}) \right)  \hfill \cr\cr
&= &\displaystyle \int_{\D^p \ts FM} \s_{\nabla}€ (\vp_0 , \ldots , 
\vp_p)^* \, \left( \pi_{1}€^{*}€({\wt \vp}^{-1\, *} (\eta)) \wedge
\wt{\vf} \right) \, . \hfill \cr
}
$$
To prove the second claim, we use 
the Stokes formula:
$$
\matrix{
&& (-1)^m \,\int_{\D^p \ts FM} \s_{\nabla}€ 
(\vp_0 , \ldots , 
\vp_p)^* (\pi_1^* (\eta) \wdg 
d \vf)  \hfill \cr\cr
&+&
\int_{\D^p \ts FM} \s_{\nabla}€ (\vp_0 , \ldots , 
\vp_p)^* (\pi_1^* (d \eta) \wdg 
\vf)  \hfill \cr\cr
&= &\displaystyle \int_{\D^p \ts FM} 
d \left(\s_{\nabla}€ (\vp_0 , \ldots , \vp_p)^* \, 
(\pi_1^* (\eta) \wdg \vf) \right)
\hfill \cr\cr
&= &\displaystyle \int_{\partial \, \D^p \ts FM} \s_{\nabla}€ (\vp_0 , \ldots , 
\vp_p)^* (\pi_1^* (\eta) \wdg 
\vf) \hfill \cr\cr
&= &\displaystyle \sum \ (-1)^i \int_{\D^{p-1} \ts FM}
\s_{\nabla}€ (\vp_0 , \ldots , 
{\check \vp}_i , \ldots , 
\vp_p)^*  (\pi_1^* (\eta) \wdg \vf) \, . \hfill \cr
}
$$
Adjusting for the sign factors, one obtains 
the stated cochain property
$$
C (d \, \vf) = (\d + \partial) \, C(\vf) \, . \qquad \hfill \Box
$$
\bigskip

To simplify the assembling of the second ingredient of our
construction, we shall take advantage of the fact that
all the cohomological information of the complex $C^* ({\Fa}_n)$
is carried by the image of the truncated Weil complex
of ${\Fg \Fl}(n, \Rb)$ via the canonical inclusion (cf. \cite{GF},
see also \cite{Bo2}, \cite{H1}). Thus,  
it suffices to work with the restriction of $C : C^* ({\Fa}_n) \ra 
TC^{*} (\Gc_M ; FM)$ to the 
subcomplex $CW^* ({\Fa}_n)$ of $C^* ({\Fa}_n)$,
generated as a graded subalgebra  by $\{ \t_j^i , R_j^i \}$, 
where $(\t_j^i )$, resp. $(R_j^i )$, is the image of the
``universal connection'' matrix,
resp. the ``universal curvature'' matrix, of the Weil complex of
${\Fg \Fl}(n, \Rb)$.

\begin{lemma} \label{Luc}
   For any torsion-free linear connection $\nabla$
on $M$, with connection form $\om_{\nabla} = (\om_j^i )$ 
and curvature form $\Om_{\nabla} = (\Om_j^i )$, one has
\begin{equation} \label{uc}
\s_{\nabla}^* (\wt{\t}_j^i) = \om_j^i \, , \qquad  
\s_{\nabla}^* (\wt{R}_j^i) = \Om_j^i \, .
\end{equation}
\end{lemma}

\noindent {\it Proof.} Since
\begin{equation} \label{curv}
R_j^i = d \, \t_j^i + \t_k^i \wedge \t_j^k \, ,
\end{equation}
the second identity is a consequence of the first. 
To prove the first, we note that by (\ref{nat}) the  
operator $\om_{\nabla} \mapsto \s_{\nabla}^* (\wt{\t})$, acting   
on the (affine) space of torsion-free 
connections on $FM$, is {\it natural}, i.e. $\Gc_M$-equivariant. 
The uniqueness result for natural operators on
torsion-free connections \cite[25.3]{K-M-S} ensures that 
the only such operator is the identity. $\qquad \Box$
\medskip

\begin{corollary} \label{polynom}
    For any   $\vf \in CW^{q} ({\Fa}_n)$ and any
    $\eta \in \Om_{c}€^{-m}€(FM) $ the integrand in formula (\ref{dchain}),
    defining the group cochain
    $\langle C_{p,m} (\vf) (\vp_0 , \ldots , \vp_p) , \, \eta \rangle $,
is a form that depends polynomially on the functions 
$\g_{jk}^i (u, \wt{\vp}_r )$ and $\g_{jk, \ell}^i (u, \wt{\vp}_s )$,
where $\, 1 \leq i, j, k, \ell \leq n$, 
$ 0 \leq r, s \leq p \, $.
\end{corollary}

\noindent{\it Proof.} By the preceding lemma applied 
to the connection (\ref{eq4.5}) together with the identity (\ref{g'}),
one has
$$
\s_{\nabla} (\vp_0 , \ldots , \vp_p)^* (\wt{\t}_j^i) = \sum_{r=0}^p \ 
t_r \, {\wt \vp}_r^* \, \om_j^i \, = \, \om_j^i \,  + \,
\sum_{r=0}^p \ t_r \, \sum_{k}€ \g_{jk}^i (u, {\wt\vp}_r) \, \t^k
\, .
$$
Taking the total exterior derivative and using (\ref{curv}), one obtains
the corresponding expression for $\s_{\nabla} (\vp_0 , \ldots , \vp_p)^* 
(\wt{R}_j^i) \, $, 
which in addition will involve $R_{j}^{i}€$ as well as 
the differentials $dt_{r}€ \, , d \, \om_j^i \, $ and
$\, d_{u}€ \, \g_{jk}^i (u, {\wt\vp}_r) \,$.

Since, by the very definition $CW^* ({\Fa}_n)$, the form ${\wt \vf}$ is
a polynomial in $ \, \wt{\t}_j^i \,$ and $\, \wt{R}_j^i \,$, 
the claim follows. $ \quad \Box$
\medskip

By composing the restriction to $\, CW^* ({\Fa}_n) \, $ of the chain map
$ C\, $ of Lemma \ref{Lchain}, (\ref{chain}) 
with the canonical map 
$$\Phi: \left( T C^{*} (\Gc_M ;  FM) , \, \d + \partial \right)
\ra \left( P C^* (\Ac) , \, b + B \right) \, ,
$$
where  $ \, \Ac : = C_c^{\ify} (FM  \semi \Gc_{M}€ ) \, $ is 
the crossed product algebra,
we obtain a new chain map 
\begin{equation} \label{newC}
\wt{C} \, := \, \Phi \circ C \, : \left( CW^* ({\Fa}_n) , \, d \right) 
\, \longra \, 
\left( P C^* (\Ac) , \, b + B \right) \, .
\end{equation}
Our next task will be to prove that the image of $\wt{C}$ 
actually lands inside the differentiable periodic cyclic
complex $P C_{\rm d}€^* (\Ac)$.

In preparation for that,
let us recall the definition of $\Phi$ \cite[III.2.$\d$]{Co}.
It involves 
the crossed product algebra $\Cc = \Bc \semi \Gc_{M}€ \, $,
where 
$$ \Bc \, = \, \Om_c^{*}€ (FM) {\hat\ot} \wedge^{*}€ \, \Cb [\Gc'_M€]  \, ,
$$
is the graded tensor product of the algebra of compactly supported
differential forms on $FM$ by the graded algebra over $\Cb$
generated by the degree $1$ anticommuting symbols 
$\d_{\vp}€ \,$, with $\vp \in \Gc_{M} \, , \vp \neq 1\, $
 and $\, \d_{1}€ = 0$. The
crossed product is taken with respect to the tensor product
action of $\Gc_{M}€$,
so that the following multiplication rules
hold:
\begin{equation}
    \matrix{
U_{\vp}^{*}€ \, \eta \, U_{\vp}€\, = \, {\wt \vp}^{*}€ \eta \, ,
\qquad \fl \, \eta \in  \Om^{*}€ (FM) \, , \hfill \cr \cr
U_{\vp}^{*}€ \, \d_{\psi}€ \, U_{\vp}€\, = \, \d_{\psi \circ \vp}€
\, - \, \d_{\vp}€ \, , \quad \fl \, \vp, \psi \in \Gc_{M}€\, . \cr
}
\end{equation}
The graded algebra $\Bc$ is endowed with the differential
$$ d \, (\eta \, \ot \d_{\vp_1} \cdots \d_{\vp_p} ) \, = \,
   d \eta \, \ot \d_{\vp_1} \cdots \d_{\vp_p} \, , \quad
   \fl \, \eta \in \Om_c^{*} \,(FM) \, ,
$$
and the crossed product algebra $\Cc$ acquires the differential
$$ d ( b \, U_{\vp}^{*} ) \, = \, (d \, b )
\, U_{\vp}^{*}€ \, - \, (-1)^{{\rm deg} \, b}€ \, b \, 
 \d_{\vp} \, U_{\vp}^{*}€ \, , \quad b \in \Bc , \, \vp \in \Gc_{M} \, .
$$
Any cochain  $\, \nu \in C^{p,m} (\Gc_M ; FM)) \,$ 
gives rise to a linear functional $\wt{\nu} \, $ on $\Cc \,$,
defined as follows, 
 $\fl \, \eta \in \Om_c^{-m} \,(FM)$,
\begin{equation} \label{fctl}
\wt{\nu} \, (\eta \, \ot \d_{\vp_1} \cdots \d_{\vp_p} 
 \, U_{\vp}^{*}) = 
\left\{ \matrix{
\langle \, \nu \, (1 , \vp_1 , 
\ldots , \vp_p) , \eta \, \rangle \quad \hbox{\rm if} \quad \vp=1 \, ,
\cr \cr
 0 \, \quad \hbox{\rm otherwise} \, .\hfill \cr
} \right. 
\end{equation}

With this notation in place, the cochain map
(\ref{newC}) is given by the formula, $\, \fl \vf \in CW^{q} ({\Fa}_n)$,
\begin{equation} \label{exp}
 \matrix{
\wt C_{p,m}€ (\vf) (f^0 \, U_{\vp^0}^* , \ldots , f^p \, U_{\vp_p}^*)
:= \hfill \cr \cr
\frac{p!}{(q+1)!} \sum_{j=0}^{q=p-m+1} \ (-1)^{j(q-j)} \, 
\wt{C_{p,m}€(\vf)} 
\left(d (f^{j+1} 
U_{\vp^{j+1}}^*) \ldots f^0 \, U_{\vp^0}^* \ldots d (f^j  U_{\vp^j}^*)
\right) .
 \cr \nonumber
}
\end{equation}

\begin{lemma} \label{coho}
    The cochain homomorphism $ \, \wt C \, = \, \Phi \circ C \,$
maps $\, CW^* ({\Fa}_n) \,$ to $\, P C€^* (\Hc_{FM}€) \, $.
\end{lemma}

\noindent {\it Proof.} We have to show
 that the cochain (\ref{exp}) is differentiable, i.e.
 of the form (\ref{dc}). 
Since $\{ \t^k , \, \om_j^i \} \,$ forms a basis of the cotangent
space at each point of $FM$, we can express the differential of
a function $f \in C_{c}€^{\ify}€ (FM)$ as
$$
df =  Y_i^j (f) \, \om_j^i +  X_k (f) \, \t^k \, 
$$
and therefore
\begin{equation} \label{dfU}
d (f \, U_{\vp}^*) =  Y_i^j (f) \, \om_j^i \, U_{\vp}^* \, + \,
 X_k (f) \, \t^k \,  U_{\vp}^* \, - \,
 f \, \d_{\vp} \, U_{\vp}^* \, .
\end{equation}
Furthermore, for the canonical form one has
\begin{equation} \label{Uth}
U_{\vp}^* \, \t^k \, U_{\vp}  = \, \t^k \, ,
\end{equation}
while for the connection form, by (\ref{g'}),
\begin{equation} \label{omU}
U_{\vp}^* \, \om_j^i \, U_{\vp} = \om_j^i + 
\g_{jk}^i (u, \wt{\vp}) \, \t^k \,  \, .
\end{equation}
Using the identities (\ref{dfU}), (\ref{Uth}), (\ref{omU}) and
Corollary \ref{polynom}, one can now achieve the proof by following
the same reasoning as in \cite[pp. 232-234]{CM2}, where
the similar result was established for the flat case.
$\quad \Box$
\medskip

The definition of the above cochain map involves the choice of 
a torsion-free connection $\nabla$. To emphasize this dependence,
we shall use the more suggestive notation
$$
{\wt C}_{\nabla}€ \, = \, \Phi \circ C_{\nabla}€ \, : \, 
CW^* ({\Fa}_n) \, \longra \, P C€^* (\Hc_{FM}€) \, .
$$  
However, the choice of the connection is
cohomologically immaterial.

\begin{lemma} \label{indep}
Let $\nb^0$, $\nb^1$ be torsion-free connections on $FM$. The
corresponding cochain homomorphisms $ \, {\wt C}_{\nabla_{0}€}€ \,$
and $ \, {\wt C}_{\nabla_{1}€}€ \,$ are cochain homotopic.
\end{lemma}

\noindent {\it Proof.} With $\s_{\wt{\nb}}  (\vp_0 , \ldots , \vp_p) 
: I \ts \D^p \ts FM \ra F^{\ify}M$ given by 
$$
\s_{\wt{\nb}} (\vp_0 , \ldots , \vp_p) (s,t,u) = \s_{(1-s) \nb^0 
(\vp_0 , \ldots , \vp_p ; t) + s 
\nb^1 (\vp_0 , \ldots , \vp_p ; t)} (u) \, ,
$$
for any $\eta \in  \Om_c^{-m} \,( FM)) \,$ we define
$$
\langle C_{\wt{\nb}} (\vf)(\vp_0 ,\ldots ,\vp_p) , \eta \rangle = 
(-1)^{\frac{m(m+1)}{2}}€
 \int_{I \ts \D^p \ts FM} 
\eta \wdg \s_{\wt{\nb}} (\vp_0 , \ldots , \vp_p)^* (\wt{\vf}) \, .
$$
Applying Stokes, one obtains
$$
\matrix{
&& \int_{I \ts 
\D^p \ts FM} \eta \wdg \, \s_{\wt{\nb}} 
(\vp_0 , \ldots , \vp_p)^* (d\vf) \hfill \cr\cr
 &+ &\displaystyle (-1)^{m}€ \, \int_{I \ts 
\D^p \ts FM} d\eta \wdg \, \s_{\wt{\nb}} 
(\vp_0 , \ldots , \vp_p)^* (\vf) \hfill \cr\cr
&= &\displaystyle \int_{I \ts \D^p \ts FM}  
d \left( \eta \wdg \s_{\wt{\nb}} (\vp_0 , \ldots , 
\vp_p)^* (\vf) \right) \, \hfill \cr\cr
&= &\displaystyle \int_{\partial I \ts \D^p \ts FM}
 \eta \wdg \s_{\wt{\nb}} (\vp_0 , 
\ldots , \vp_p)^* (\vf) \hfill \cr\cr
&+ &\displaystyle \int_{I \ts \partial \D^p \ts FM} 
\eta \wdg \s_{\wt{\nb}} (\vp_0 , \ldots , 
\vp_p)^* (\vf) \, . \hfill \cr
}
$$
After adjusting for the sign factors, this gives the desired
homotopy formula
$$
C_{\nb^1} - C_{\nb^0} = C_{\wt{\nb}} \circ d - (\partial + \d) 
\circ C_{\wt{\nb}} \, . \qquad \Box
$$
\medskip

Let $K \, \sbs \, GL(n, \Rb) \, $ be a compact subgroup. By restricting
the map $C_{\nabla}€$
to the subcomplex $ CW^* ({\Fa}_n , K) \, $ of $K$-basic elements
in $\,CW^* ({\Fa}_n) \,$, one obtains a chain map
$ \, C^{K}_{\nabla}€ : C^* ({\Fa}_n , K) \ra C^{*} 
(\Gc_M ; FM/K) \, $. Similarly, the map $\Phi$ restricts to
$\, \Ph^{K}€ : T C^{*} (\Gc_M ;  FM/K) 
\ra  P C^* (\Ac^{K}€) \, $, where we identify
$$ \, \Ac^{K}€ \simeq \Ac_{FM/K} := C_c^{\ify} (FM/K  \semi \Gc_{M}€ ) \, .
$$
By composition, one gets the relative chain map  
\begin{equation} \label{relC}
\wt{C}^{K}_{\nabla}€ \, := \, 
\Phi^{K}€ \circ C_{\nabla}€ \, : \, CW^* ({\Fa}_n , K)  \longra 
 P C^* (\Hc_{FM} , K) \, .
\end{equation}

We are now ready to conclude the proof of our main result, which
is a refinement of Theorem \ref{MI} and provides a geometric 
realization of the Gelfand-Fuchs characteristic classes in
the framework of cyclic cohomology.

\begin{theorem} \label{vE} Let $M$ be a smooth $n$-dimensional
    manifold endowed with a torsion-free linear connection $\nabla$
    and let $K \, \sbs \, O(n) \, $ be a compact subgroup.
    The cochain homomorphism (\ref{relC})
 induces isomorphism in cohomology,
 \begin{equation} \label{gf} 
 \g^{*}_{K}€ : \sum_{i \equiv * \, {\rm mod} \, 2} \, 
 H^{i}€ ({\Fa}_n , K) \build \longra_{}^{\simeq} 
 H C_{\rm per}^* (\Hc_{FM} , K)\, .
\end{equation}
\end{theorem}

\noindent {\it Proof.} Returning to the setting of the proof of
Theorem \ref{MI}, we choose an open cover of $M$ 
by domains of local charts $\Uc = \{ V_i \}_{1 \leq i \leq r} $,
together with a partition of unity $\{  \chi_{i}^{2} \}_{1 \leq i \leq 
r} \,$ subordinate to the cover $\Uc \,$. We then form the
flat manifold 
$N \, = \, \coprod_{i=1}^r \, V_{i} \ts \{ i \} \, $
and the corresponding smooth \'etale groupoid
$\Gc_{\Uc}$, such that the algebras
$\, \Ac =  C_c^{\ify} (FM \semi \Gc_M ) \,$ and
$\, \Bc = C_c^{\ify} (FN \semi \Gc_{\Uc}) \,$ are Morita equivalent.

It suffices to show that by composing
$\, {\wt C}_{\nabla}€ \,$ with the cochain equivalence
$\Ps^{*}€ : C_{d}^{*}€(\Ac) \ra C_{d}^{*}€(\Bc) \,$ one gets
a quasi-isomorphism. In turn, the composite map
$\, \Ps^{*}€ \circ {\wt C}_{\nabla}€ \,$ can be seen to be
homotopic to
the map 
$$\, {\wt C}_{\wt \nabla}€  : \,CW^* ({\Fa}_n) \,
\longra \, P C€^* (\Hc_{FN}€) \, ,
$$ 
corresponding to the
connection $\, {\wt \nabla} \,$ on $N$ obtained by restricting
$\nabla$ to each $V_i \, , 1 \leq i \leq r$. We are thus
reduced to proving the statement on $N$, where
one can apply Lemma \ref{indep} and replace $ {\wt \nabla} $ by 
the trivial flat connection $\nabla_0$. In that case the statement
follows from \cite[\S 7, Theorem 11]{CM2}, after noticing that
the cochain map of \cite[Lemma 8]{CM2} and its variant
$\, {\wt C}_{\nabla€_0} \,$ employed here
are obviously homotopic. 

The relative case can be proven by restriction to $K$-invariants.
$\qquad \Box$
\medskip
 
By composing the isomorphism (\ref{gf}) with
the natural forgetful homomorphism, one obtains a characteristic
homomorphism
\begin{equation} \label{char} \chi_{O(n)}^* \, :
H^{*}€ ({ \frak a}_n , O(n)) \build \longra_{}^{\g_{O(n)}^{*}€}
  H C_{\rm per}^* (\Hc_{FM} , O(n)) \ra 
    H C_{\rm per}^{*} (\Ac_{PM}€)_{(1)}€  \,  ,
\end{equation}
where $PM \, = \, FM/O(n)$; it lands in the cyclic cohomology
component $H C_{\rm per}^{*} (\Ac_{PM}€)_{(1)}€ \,$
corresponding to the localization at the identity. One can
further compose $ \chi_{O(n)}^* \,$ with the restriction homomorphism 
$$ \i_{M}^{*}€ : H C_{\rm per}^{*} (\Ac_{PM}€)_{(1)}€
\longra H C_{\rm per}^{*} (C_{c}^{\ify}€(PM)) \, ,
$$
corresponding to the natural inclusion
$ \i_{M}€ : C_{c}^{\ify}€(PM) \ra \Ac_{PM}€ \, .$
After identifying the target to the twisted cohomology
(by the orientation sheaf) 
$$\, H_{\tau}^{*}€ (PM) \simeq H_{\tau}^{*}€ (M) \, , 
$$
one gets a natural homomorphism
$$ \chi_{M}^* : H^{*}€ ({ \frak a}_n , O(n)) \longra 
H_{\tau}^{*}€ (M) \, .
$$

\begin{proposition} The homomorphism $ \chi_{M}^* : H^{*}€ ({ \frak a}_n , O(n)) 
    \longra H_{\tau}^{*}€ (M) \, $ is the classical
 characteristic homomorphism, expressing the Pontryagin classes
 $ (p_{i_{1}€}€ \cdots  p_{i_{k}€}) \, (M) \,$  of $M$ as 
 images of the universal Chern classes $ c_{2i_{1}€}€ \cdots  c_{2i_{k}€}€ 
  \in H^{*}€ ({ \frak a}_n , O(n)) \, ,
  \quad 2i_{1}€ + \ldots + 2i_{k}€ \leq n $.
\end{proposition}

\noindent {\it Proof.} Let $\nabla$ 
be a torsion-free connection. 
For any $\vf \in CW^{q}€ ({\Fa}_n , O(n)) \,$ and any
$ f_0 , f_{1} , \ldots , f_m  \in  C_{c}^{\ify}€(FM/O(n))
\,$, one has (up to sign)
$$
\langle C_{0, -m}(\vf)\, (1) \, , \,
f_0 \, df_{1}€ \wedge \ldots \wedge df_m \rangle = 
 \int_{ PM} \, f_0 \, df_{1}€ \wedge \ldots \wedge df_m 
\wdg \s_{\nb}^{*}€ (\wt{\vf}) \, .
$$
Thus, the statement is an immediate consequence of
Lemma \ref{Luc}. 
$ \qquad \Box$    

\section{Application to the transverse index formula}

We shall now indicate how the preceding results apply to 
the cohomological index formula for the transverse
fundamental class in $K$-homology (cf. \cite{CM1}, \cite{CM2},
\cite{CM4}), in the setting of diffeomorphism invariant
geometry.

With the manifold $M$ assumed to be oriented, we first recall the
definition of the spectral triple that encodes its 
diffeomorphism invariant fundamental class.
One starts by forming 
the bundle of local metrics (cf. \cite{C4}),  
$ \pi: PM \, = \,  F^{+}€M / SO (n) \, \ra \, M $, where
 $F^{+}€M$ is bundle of oriented frames on $M$.
The vertical subbundle $\Vc \sbs TP$, $\Vc =\Ker \pi_*$, 
carries natural inner products on each of its
fibers, determined solely by the choice of a $GL^+(n, \Rb)$--invariant
Riemannian metric on 
the symmetric space $GL^{+}€(n,  \Rb) / SO (n)$. At the same time,
the quotient bundle $\Nc = (TP)/ \Vc$ comes equipped with its own,
tautologically defined, Riemannian structure: every $p\in P$ is an
Euclidean structure on $T_{\pi (p)} (M)$ which is identified to 
$\Nc_{p}€$ via $\pi_*$. 
By the naturality of the above construction, 
the pseudogroup of orientation 
preserving local diffeomorphisms
$\G^{+}_{M}€$, acting by prolongation on $PM$, 
preserves the ``para-Riemannian'' structure thus
defined. The algebra of ``coordinates'' of the spectral triple
is the convolution algebra of 
the smooth \'etale groupoid $ PM \bar{\semi} \G^{+}_{M}€ $,
$$ 
     {\Ac}_{PM}€ = C_c^{\ify} ( PM \bar{\semi} \G^{+}_{M}€)   \, .
$$
The Hilbert space of the spectral triple is
$$\, L^{2}€({\wedge}^{\cdot}€ \Vc^{*}€  \ot 
 {\wedge}^{\cdot}€ \Nc^{*}€  , \,  vol_{P}€) \, ,
 $$
where $vol_{P}$ denotes the canonical $\G^{+}_{M}€$-invariant
volume form on $P$. The algebra  ${\Ac}_{PM}€$ acts on this
space by multiplication operators
$$
    ( (f \, U_{\psi}^*) \, \xi) (p) = f(p) \, \xi (\wt{\psi} (p)) \quad \fl
\, p
\in PM \, , \ \xi \in L^{2}€(P) \, , \ f \, U_{\psi}^* 
\in {\Ac_{PM}€} \, . 
$$
To complete the description of the spectral triple,
we need to define the operator $D$, representing the
$K$-homology orientation class of $M$
in a diffeomorphism invariant fashion.
It is given by the identity $\, Q = D |D| \,$,
where the {\it hypoelliptic signature operator} 
$Q$ is defined as a graded sum
\begin{equation} \label{Q}
    Q = (d_V^* \, d_V - d_V \, d_V^*) \op \g \, (d_H + d_H^*) \, ;
\end{equation}
$d_V$ denotes the vertical exterior derivative and
$d_H $ stands for the horizontal exterior differentiation with
respect to a fixed torsion-free connection $\nabla$.
When $n \equiv 1 \, {\rm or} \, 2 \, ({\rm mod} \,  4)$, 
in order for the vertical component to make sense,
one has to replace $PM$ by $ PM \times S^1$ 
so that the dimension of the vertical fiber
stays even. 
As shown in \cite{CM1}, all the
commutators $\, [D,a] \,$, with $\, a \in \Ac_{PM}€ \, ,$ 
are bounded. Furthermore, for any
$f \in C_{c}€^{\infty}€(P) \, $ and any $  \lb \notin \Rb $, the
local resolvent
$\, f(D-\lb)^{-1} \, $ is $p$-summable, for every $p$ that exceeds 
 the Hausdorff
dimension $\, d = {n(n+1) \over 2} + 2n \,$
of $PM$ viewed as a Cartan-Carath\'eodory metric space.

As a $K$-homology class, the operator $D$ 
determines an additive map from the $K$-theory group
$ K_{*}€ (\Ac_{PM}€)$ to  $\Zb $, via the
familiar index pairing:
\begin{itemize}
    \item[(0)] in the {\it graded} (or {\it even}) case,
    $$ {\rm Index}_{D}€ ([e]) = {\rm Index} \, (e D^{+}€e) \,  , 
    \quad e^{2}€ = e \in \Ac_{PM}€ \, ;
    $$
    \item[(1)] in the {\it ungraded} (or {\it odd}) case,
    $$ {\rm Index}_{D}€ ([u]) = {\rm Index}\, (P^{+}€ u P^{+}€) \,  , \quad
    u \in GL_{1}€(\Ac_{PM}€) \, ,
    $$
    where $P^{+}€ = {1+F \over 2}$, with  $F = \Sign \, D$ .
\end{itemize} 
In cohomological terms, the index pairing can be expressed as
as a pairing between cyclic (co)homological classes :
\begin{equation} \label{index}
 {\rm Index}_{D}€ (\kappa) = \lgl {ch}_* (D) , \, {\rm ch}^* (\kappa) \rgl 
\qquad \fl \, \kappa \in  K_{*}€ (\Ac_{PM}€)  .
\end{equation}
While $ {\rm ch}^* (\kappa) \in  H C^{\rm per}_{*} (\Ac_{PM}€) \, $ 
is easy to express in terms of the ``difference idempotent'' determined
by the $K$-theory class $\, \kappa \, $, 
the Chern character class in $K$-homology 
 $ \, {ch}_* (D) \in H C_{\rm per}^{*} (\Ac_{PM}€) \,$
 has a more involved definition (cf. \cite{C1}). In 
the odd case it is given by the cyclic cocycle
\begin{equation}  \label{Kch}
\tau_{F}€ (a^0 ,\ldots ,a^p) = \Trace \, (a^0 [F,a^1] \ldots [F,a^p]) 
\, , \qquad  a^j \in \Ac_{PM}€ \, ,
\end{equation}
where $p$ is any odd integer
exceeding the dimension  $\, d = {n(n+1) \over 2} + 2n \,$
of the spectral triple;
in the even case the trace gets replaced by the graded trace 
$\Trace_{s}€$ and $p$ is even.
As such, 
the cocycle (\ref{Kch}) is inherently difficult to compute, because
it involves the operator trace. 

In \cite[Part I]{CM1} we used 
the hypoelliptic calculus on Heisenberg manifolds
adapted to the para-Riemannian structure 
of the manifold $PM$ to prove that
the spectral triple constructed above 
fulfills the hypotheses of the
operator theoretic local index theorem of \cite[Part II]{CM1}. 
In particular, a pseudodifferential operator $T$ in the
aforementioned calculus admits a
locally computable residue of Wodzicki-Guillemin-Manin-type
\begin{equation} \label{res}
{\int \!\!\!\!\!\! -} T = {1 \over (2 \pi)^{d}€} 
\int_{PM}€ \, {\rm res}_{T}€(u) \, vol_{PM}€(u)
\end{equation}
where $\, {\rm res}_{T} (u) $ is the function obtained by integrating
the symbol of $T$ of critical order $ - d \,$, 
$\, \s^{T}_{- d}€ (u, \xi) \,$,
against the Liouville measure of the sphere $\, ||\xi||^{'}€ = 1 $
corresponding to the intrinsic quartic metric of $PM$.
Specializing the local index formula of \cite[Part II]{CM1} to
the above spectral triple, assumed for definiteness to be 
odd-dimensional, one obtains that the
 Chern character $ {ch}_* (D) \in H C_{\rm per}^{*} (\Ac_{PM}€) \,$
can be represented by 
 the cocycle $ \Phi_{Q}€ = \{ \ph_q  \}_{q=1,3,\ldots} \,$ 
 defined as follows:
\begin{eqnarray} \label{phi}
    \matrix{
\ph_q (a^0 ,\ldots ,a^q) =  \hfill  \cr \cr
 \sum_k (-1)^{\vert k \vert} \, \sqrt{2i} (k_1 ! \cdots k_q !)^{-1} 
((k_1 +1) \ldots (k_1 + \ldots + k_q +q))^{-1}  
\G ( \vert k \vert + {q \over 2} ) \times  
}
\end{eqnarray}
$$
{\int \!\!\!\!\!\! -}  a^0 [Q, a^1]^{(k_1)} \ldots
[Q, a^q]^{(k_q)}  \, \vert Q \vert^{-q-2\vert k\vert} \, ,
\quad a^j \in \Ac_{PM}€ \, ,
$$
where  we used the abbreviations $\vert k \vert = k_1 +\ldots + k_q$
and
$$T^{(i)} = \nb^i (T) \qquad {\rm and} \qquad \nb (T) = D^2 T
- TD^2 \, .
$$
Only finitely many terms in the above sum may not vanish and only
finitely many components of $\Phi_{Q}€$ are nonzero.
\medskip

While the expression (\ref{phi}) is algorithmically computable in
principle,
its actual computation is prohibitively difficult to perform in practice.
However, as we showed in \cite{CM2} for the case of the flat connection, 
from the cohomological
standpoint the answer is as reasonable as it could possibly be.
(See also \cite{Pe} for a relevant illustration.) 

We are now in a position to remove the flatness assumption
from the statement of our transverse index theorem 
\cite[\S 9, Theorem 5]{CM2} and formulate it in full generality,
for a hypoelliptic signature operator formed with 
an arbitrary torsion-free coupling connection.
\smallskip

\begin{theorem} Let $Q = D |D| $ be the hypoelliptic signature
operator on $PM$ associated to a torsion-free connection
$\nabla$. The identity
component of its Chern character, $ {ch}_* (D)_{(1)}€ 
\in H C_{\rm per}^{*} (\Ac_{PM}€)_{(1)}€ \,$, is
the image under the characteristic homomorphism (\ref{char})
of a universal class 
$\Lc_{n} \, \in \, H^{*}€ ({ \frak a}_n , SO(n)) \,$,
$$ {ch}_* (D)_{(1)}€ \, = \, \chi_{SO(n)}^* \, (\Lc_{n} ) \, .
$$
In particular, the class
$ {ch}_* (D)_{(1)}€$ can be represented by a cocycle built out of the 
connection form $\om_{\nabla}€ = ( \om^i_j )$, 
its curvature form $\Om_{\nabla}€ = ( \Om^i_j )$
and the corresponding displacement
functions on the jet groupoid of $ FM \bar{\semi}\G_{M}€$, 
$ \, \g_{jk}^i \,$ and $\, \g_{jk , \, \ell}^i \, \equiv
 \, X_{\ell}€\, \g_{jk}^i \,$ ,
$ 1 \leq i, j, k, \ell \leq n \,$.

\end{theorem}

\noindent {\it Proof.} In view of the 
formula (\ref{phi}) for $ch_{*}€(D )$, it suffices to show that 
any cochain on ${\Ac}_{PM}$ of the form,
$$ 
\ph (a^0 , \ldots , a^q) = \, \int \!\!\!\!\!\! - \, a^0 [Q, a^1]^{(k_1)}
\ldots
[Q , a^q]^{(k_q)} \, |Q|^{- (q + \vert 2k \vert)} \,  ,
$$
with $\, a^j \, = \, f^j \, U_{\psi_{j}€}^{*}€ \,
\in \, {\Ac}_{PM} \, , \quad j = 1, \ldots , q \, , $ such that 
$ \, \psi_{q}€ \circ \cdots \circ \psi_{1}€ \circ \psi_{0}€ = 1 \, $,
belongs to the range of the characteristic map $\chi_{SO(n)}^*€ $.
Using the results in \cite[\S 9]{CM2}, specifically Lemma 1 and its
corollary, one can write such a cochain as a sum of 
cochains of the form 
$$
 \int \!\!\!\!\!\! - \, a^0 \, h^{1} (a^1) \cdots h^{q} (a^q) \, R \, 
, \quad h^{1}€ , \ldots , h^{q}€  \in \Hc_{FM}€ \, ,
$$
with $R$ a pseudodifferential operator in the hypoelliptic
calculus. By (\ref{res}), each such expression is indeed
of the desired form
$$
 \tau \left( a^0 \, h^{1} (a^1) \cdots {\wt h}^{q} (a^q) \right) \, ,
$$
once the local residue function $ {\rm res}_{R} $ is absorbed into
$ {\wt h}^{q} \, = \, \b ({\rm res}_{R}€) \, h^{q} $.  

Finally, the second assertion follows from  Corrolary \ref{polynom} 
and Theorem \ref{vE}. $ \qquad \Box $

\end{document}